\documentclass[12pt,letterpaper]{article}
\usepackage{amsthm,amssymb}
\usepackage{amsfonts}
\makeatletter
\usepackage[dvips]{color}
\usepackage{colordvi,multicol}

\newtheorem{thm}{\textbf Theorem}[section]
\newtheorem{lem}[thm]{\textbf Lemma}
\newtheorem{rem}[thm]{\textbf Remark}

\newtheorem{prop}[thm]{\textbf Proposition}
\newtheorem{defin}[thm]{\textbf Definition}
\newtheorem{assum}[thm]{\textbf Assumption}

\newtheorem{exa}[thm]{\textbf Example}
\newcommand{\be}{\begin{eqnarray*}}
\newcommand{\ee}{\end{eqnarray*}}

\begin{document}

\title{\bf FUKUSHIMA'S DECOMPOSITION FOR DIFFUSIONS ASSOCIATED WITH SEMI-DIRICHLET FORMS}
  \author{LI MA\\
Department of Mathematics\\
Hainan Normal University\\
Haikou, 571158, China\\
mary\hskip -0.09cm\_\hskip 0.05cm henan@yahoo.com.cn\\
\\
ZHI-MING MA\\
Institute of Applied Mathematics\\
AMSS, Chinese Academy of Sciences\\
Beijing, 100190, China\\
mazm@amt.ac.cn\\
\\
WEI SUN
            \\ Department of Mathematics and Statistics\\
             Concordia University\\
             Montreal, H3G 1M8, Canada\\
             wsun@mathstat.concordia.ca}
   \date{}
\maketitle

\begin{abstract}
\noindent Diffusion processes associated with semi-Dirichlet forms are studied in the paper. The main results are Fukushima's decomposition for the diffusions and a transformation formula for the corresponding martingale part of the decomposition. The results are applied to some concrete examples.
\end{abstract}
%%%%%%%%%%%%%%%%%%%%%%%%%%%%%%%%%%%%%%%%%%%%%%%%%%%%%%%%%%%%%%%%%%%%%%%
\section[short title]{Introduction}

It is well known that Doob-Meyer decomposition and It$\hat{\rm o}$'s
formula are essential in the study of stochastic dynamics. In the
framework of Dirichlet forms, the celebrated Fukushima's
decomposition and the corresponding transformation formula play the
roles of Doob-Meyer decomposition and It$\hat{\rm o}$'s formula,
which are available for a large class of processes that are not semi-martingales. The classical
decomposition of Fukushima was originally established for regular
symmetric Dirichlet forms (cf. \cite{Fu79} and \cite[Theorem
5.2.2]{Fu94}). Later it was extended to the non-symmetric and
quasi-regular cases, respectively (cf. \cite[Theorem 5.1.3]{oshima}
and \cite[Theorem VI.2.5]{MR92}). Suppose that $({\cal E},D({\cal
E}))$ is a quasi-regular Dirichlet form on $L^2(E;m)$ with
associated Markov process $((X_t)_{t\ge 0}, (P_x)_{x\in
E_{\Delta}})$ (we refer the reader to \cite{MR92}, \cite{Fu94} and
\cite{MR95} for notations and terminologies of this paper). If $u\in
D({\cal E})$, then there exist unique martingale additive functional
(MAF in short) $M^{[u]}$ of finite energy and continuous
additive functional (CAF in short) $N^{[u]}$ of zero energy
such that
$$
\tilde{u}(X_t)-\tilde{u}(X_0)=M^{[u]}_t+N^{[u]}_t,
$$
where $\tilde{u}$ is an ${\cal E}$-quasi-continuous $m$-version of
$u$ and the energy of an AF $A:=(A_t)_{t\ge 0}$ is defined to be
\begin{equation}\label{1}
e(A):=\lim_{t\rightarrow 0}\frac{1}{2t}E_m[A^2_t]
\end{equation}
whenever the limit exists in $[0,\infty]$.

The aim of this paper is to establish Fukushima's decomposition for
some Markov processes associated with semi-Dirichlet forms. Note
that the assumption of the existence of dual Markov process plays a
crucial role in all the Fukushima-type decompositions known up to
now. In fact, without that assumption, the usual definition
(\ref{1}) of energy of AFs is questionable. To tackle this
difficulty, we employ the notion of local AFs (cf. Definition
\ref{localAF} below) introduced in \cite{Fu94} and
introduce a localization method to obtain Fukshima's decomposition for
a class of diffusions associated with semi-Dirichlet forms.  Roughly
speaking, we prove that for any $u\in {D(\mathcal{E})}_{loc}$, there exists a unique decomposition
$$
\tilde{u}(X_{t})-\tilde{u}(X_{0})=M^{[u]}_{t}+N^{[u]}_{t},
$$
where $M^{[u]}\in{\mathcal{M}}^{[\![0, \zeta[\![}_{loc}$ and
$N^{[u]}\in {\mathcal{N}}_{c,loc}$. See Theorem \ref{thm3.2} below for the involved
notations and a rigorous statement of the above decomposition.

 Next,  we develop a transformation formula
 of local MAFs.  Here we encounter the difficulty that there is no
LeJan's transformation rule available for semi-Dirichlet forms. Also
we cannot replace a
$\gamma$-co-excessive
 function  $g$ in the Revuz correspondence (cf. (\ref{e:revuz}) below) by an
  arbitrary $g\in\mathcal{B}^{+}(E)\cap
D(\mathcal{E})$, provided the corresponding smooth measure is not of
finite energy integral. Borrowing some ideas of \cite[Theorem 5.4]{Kim} and
\cite[Theorem 5.3.2]{oshima}, but putting more extra efforts, we are
able to build up an analog of LeJan's formula (cf. Theorem
\ref{thm4.6} below). By virtue of LeJan's formula developed in
Theorem \ref{thm4.6}, employing again the localization method
developed in this paper, finally we obtain a transformation formula
of local MAFs for semi-Dirichlet forms in Theorem \ref{Ito}.

The main results derived in this paper rely heavily on the potential
theory of
semi-Dirichlet forms. Although they are more or less parallel to
those of symmetric Dirichlet forms, we cannot find explicit
statements in literature. For the solidity of our results, also
for the interests by their own, we checked and derived in detail
some results on potential theory and positive continuous AFs (PCAFs in short) for semi-Dirichlet forms.
These results are presented in Section \ref{Appendix}  at the end of
this paper as an Appendix. In particular, we would like to draw the attention of the
readers to two new results, Theorem
\ref{thm34} and Lemma \ref{lnm}.

The rest of the paper is organized as follows. In Section
\ref{Sec:Fukushima}, we derive Fukushima's decomposition. Section
\ref{sec:transform} is devoted to the transformation formula. In Section \ref{Sec:example}, we apply our main results
to some concrete examples.  In these examples the usual Doob-Meyer
decomposition and It$\hat{\rm o}$'s formula for semi-martingales are
not available. Nevertheless we can use our results to preform
Fukushima's decomposition and apply the transformation formula in the
semi-Dirichlet forms setting. Section \ref{Appendix} is the Appendix
consisting of some results on potential theory and PCAFs in
the semi-Dirichlet forms setting.

%%%%%%%%%%%%%%%%%%%%%%%%%%%%%%%%%%%%%%%%%%%%%%%%%%%%%%%%%%%%%%%%%%%%%%%%%%%%%%%%%%%%%%%%%%%%%%%%%%
\section [short title]{Fukushima's decomposition}\label{Sec:Fukushima}
\setcounter{equation}{0}  We consider a quasi-regular semi-Dirichlet
form  $({\cal E},D({\cal E}))$ on $L^2(E;m)$, where $E$ is a
metrizable Lusin space (i.e., topologically isomorphic to a Borel subset of a
complete separable metric space) and $m$ is a $\sigma$-finite positive
measure on its Borel $\sigma$-algebra ${\cal B}(E)$. Denote by
$(T_{t})_{t\geq0}$ and $(G_{\alpha})_{\alpha\geq0}$ (resp.
$(\hat{T}_{t})_{t\geq0}$ and $(\hat{G}_{\alpha})_{\alpha\geq0}$) the
semigroup and resolvent (resp. co-semigroup and co-resolvent)
associated with $({\cal E},D({\cal E}))$. Let ${\bf M}=(\Omega,{\cal
F},({\cal F}_t)_{t\ge 0}, (X_t)_{t\ge 0},(P_x)_{x\in E_{\Delta}})$
be an $m$-tight special standard process which is properly
associated with $({\cal E},D({\cal E}))$ in the sense that $P_tf$ is
an ${\cal E}$-quasi-continuous $m$-version of $T_tf$ for all $f\in
{\cal B}_b(E)\cap L^2(E;m)$ and all $t>0$, where $(P_t)_{t\ge 0}$
denotes the semigroup associated with ${\bf M}$ (cf. \cite[Theorem
3.8]{MR95}). Below for notations and terminologies related to
quasi-regular semi-Dirichlet forms we refer to \cite{MR95} and
Section \ref{Appendix} of this paper.

 Recall that a positive measure $\mu$ on $(E,
{\cal B}(E))$ is called {\it smooth} (w.r.t. $({\cal E},D({\cal
E}))$), denoted by $\mu \in S,$  if $\mu(N)=0$ for each ${\cal
E}$-exceptional set $N\in {\cal B}(E)$ and there exists an ${\cal
E}$-nest $\{F_k\}$ of compact subsets of $E$ such that
$$
\mu(F_k)<\infty\ {\rm for\ all}\ k\in\mathbb{N}.
$$
A family $(A_t)_{t\ge 0}$ of functions on $\Omega$ is called an {\it
additive functional} (AF in short) of ${\bf M}$ if:

(i) $A_t$ is ${\cal F}_t$-measurable for all $t\ge 0$.

(ii) There exists a defining set $\Lambda\in{\cal F}$ and an
exceptional set $N\subset E$ which is ${\cal E}$-exceptional such
that $P_x[\Lambda]=1$ for all $x\in E\backslash N$,
$\theta_t(\Lambda)\subset\Lambda$ for all $t>0$ and for each
$\omega\in \Lambda$, $t\rightarrow A_t(\omega)$ is right continuous
on $(0,\infty)$ and has left limits on $(0,\zeta(\omega))$,
$A_0(\omega)=0$, $|A_t(\omega)|<\infty$ for $t<\zeta(\omega)$,
$A_t(\omega)=A_{\zeta}(\omega)$ for $t\ge\zeta(\omega)$, and
\begin{eqnarray}\label{additivity}
A_{t+s}(\omega)=A_t(\omega)+A_s(\theta_t\omega),~~~~\forall~ s,t\ge
0.
\end{eqnarray}
Two AFs $A=(A_t)_{t\ge 0}$ and $B=(B_t)_{t\ge 0}$ are said to be
equivalent, denoted by $A=B,$ if they have a common defining set
$\Lambda$ and a common exceptional set $N$ such that
$A_t(\omega)=B_t(\omega)$ for all $\omega\in\Lambda$ and $t\ge 0$.
An AF $A=(A_t)_{t\ge 0}$ is called a continuous AF (CAF in short) if
$t\rightarrow A_t(\omega)$ is continuous on $(0,\infty).$ It is
called a positive continuous AF (PCAF in short) if $A_t(\omega)\ge
0$ for all $t\ge 0$, $\omega\in\Lambda$.

In the theory of Dirichlet forms, it is well known that there is a one to
one correspondence between the family of
all equivalent classes of PCAFs and the family $S$ (cf. \cite{Fu94}). In \cite{F01},
Fitzsimmons extended the smooth measure characterization of PCAFs
from the Dirichlet forms setting to the semi-Dirichlet forms setting. Applying \cite[Proposition 4.12]{F01}, following the arguments of
\cite[Theorems 5.1.3 and 5.1.4]{Fu94} (with slight modifications by
virtue of \cite{MR95,MR1995,kuwae} and \cite[Theorem 3.4]{f}), we
can also obtain a one to one correspondence between the family of
all equivalent classes of PCAFs and the family $S$. The correspondence, which is referred to as {\it Revuz
correspondence}, is described in the following lemma.
\begin{lem}\label{lem:Revuz}
Let $A$ be a PCAF. Then there exists a unique $\mu\in S$,
which is referred to as the {\it Revuz measure} of $A$ and is
denoted by $\mu_A,$ such that:

For any $\gamma$-co-excessive function $g$ $(\gamma\geq0)$ in
$D({\cal E})$ and $f\in \mathcal{B}^{+}(E)$,
\begin{equation}\label{e:revuz}
 \lim_{t\downarrow0}\frac{1}{t}E_{g\cdot
 m}((fA)_{t})=<f\cdot\mu,\tilde{g}>.
\end{equation}
Conversely, let $\mu\in S$, then there exists a
unique (up to the equivalence) PCAF $A$ such that $\mu =\mu_A.$
\end{lem}
See Theorem \ref{thm2.25} in the Appendix at the end of this paper for
more descriptions of the Revuz correspondence (\ref{e:revuz}).

From now on we suppose that $(\mathcal{E},D(\mathcal{E}))$ is a
quasi-regular local semi-Dirichlet form on $L^{2}(E;m)$. Here
``local" means that ${\cal E}(u,v)=0$ for all $u,v\in D({\cal E})$
with ${\rm supp}[u]\cap{\rm supp}[v]=\emptyset$.  Then,
$(\mathcal{E},D(\mathcal{E}))$ is properly associated with a
diffusion process  ${\bf M}=(\Omega,{\cal F},({\cal F}_t)_{t\ge 0},
(X_t)_{t\ge 0},(P_x)_{x\in E_{\Delta}})$ (cf. \cite[Theorem
4.5]{kuwae}). Here ``diffusion" means that ${\bf M}$ is a right
process satisfying
$$
P_x[t\rightarrow X_t\ {\rm is\ continuous\ on\ } [0,\zeta)]=1\
{\rm for\ all}\ x\in E.
$$
Throughout this paper, we fix a function $\phi\in L^{2}(E;m)$ with
 $0<\phi\leq1$ $m{\textrm{-}}a.e.$ and set $h=G_{1}\phi$,
 $\hat{h}=\hat{G}_{1}\phi$.  Denote $\tau_B:=\inf\{t>0\,|\, X_t\notin B\}$ for $B\subset E$.

Let $V$ be a quasi-open subset of $E$. We denote by
$X^V=(X^V_t)_{t\ge 0}$ the part process of $X$ on $V$ and denote by
$(\mathcal{E}^{{V}},D(\mathcal{E})_{{V}})$ the part form of
$(\mathcal{E},D(\mathcal{E}))$ on $L^2(V;m)$. It is known that $X^V$
is a diffusion process and
$(\mathcal{E}^{{V}},D(\mathcal{E})_{{V}})$ is a quasi-regular local
semi-Dirichlet form (cf. \cite{kuwae}). Denote by
$({T}_{t}^{V})_{t\ge 0}$, $(\hat{T}_{t}^{V})_{t\ge 0}$,
$({G}_{\alpha}^{V})_{\alpha\ge 0}$ and
$(\hat{G}_{\alpha}^{V})_{\alpha\ge 0}$ the semigroup, co-semigroup,
resolvent and co-resolvent associated with
$(\mathcal{E}^{{V}},D(\mathcal{E})_{{V}})$, respectively. One can
easily check that $\hat{h}|_{V}$ is 1-co-excessive w.r.t.
$(\mathcal{E}^{{V}},D(\mathcal{E})_{{V}})$. Define
$\bar{h}^V:=\hat{h}|_{V}\wedge\hat{G}^{{V}}_{1}\phi$. Then
$\bar{h}^V\in D(\mathcal{E})_{{V}}$ and $\bar{h}^V$ is
 1-co-excessive.

 For an AF $A=(A_t)_{t\ge 0}$ of $X^V$, we define
\begin{equation}\label{ebn}
e^V(A):=\lim_{t\downarrow0}{1\over{2t}}E_{\bar{h}^V\cdot
m}(A_{t}^2)
\end{equation}
whenever the limit exists in $[0,\infty]$.
Define
 \begin{eqnarray*}
\dot{\mathcal{M}}^V
 &:=&\{M\,|\, M\ \mbox{is an AF of}\ X^{V},\ E_x(M^{2}_t)<\infty, E_x(M_t)=0\\
 &&\ \ \ \ \mbox{for}\ {\rm all}\ t\ge0\ {\rm and}\ {\cal
E}{\textrm{-}q.e.}\ x\in V,
e^V(M)<\infty\},
 \end{eqnarray*}
  \begin{eqnarray*}
 \mathcal{N}^V_{c}
 &:=&\{N\,|\, N\ \mbox{is a CAF of}\ X^{V},E_x(|N_t|)<\infty\ \mbox{for}\ {\rm all}\ t\ge0\\
   &&\ \ \ \ {\rm and}\ {\cal
E}{\textrm{-}q.e.}\ x\in V, e^V(N)=0\},
   \end{eqnarray*}
\begin{eqnarray*}
\Theta&:=&\{\{{V_n}\}\,|\, {V_n}\ \mbox{is}\ {\cal E}{\textrm{-}}\mbox{quasi}
{\textrm{-}}\mbox{open},\ {V_n}\subset V_{n+1}\ {\cal E}
{\textrm{-}q.e.}, \\
         &&\ \ \ \ \ \ \ \ \ \ \ \forall~n\in \mathbb{N},\
         \mbox{and}\ E=\cup_{n=1}^{\infty}{V_n}\ {\cal E}{\textrm{-}q.e.}\},
\end{eqnarray*}
and
\begin{eqnarray*}
{D(\mathcal{E})}_{loc}
&:=&\{u\,|\,\exists\ \{V_n\}\in\Theta\ \mbox{and}\ \{u_n\}\subset D(\mathcal{E})\nonumber\\
  &&\ \ \ \ \ \ \ \mbox{such that }\ u=u_n\ m{\textrm{-}a.e.}\ \mbox{on}\ V_n,
  ~\forall~n\in \mathbb{N}\}.
\end{eqnarray*}
For our purpose we shall employ the notion of local AFs introduced
in \cite{Fu94} as follows.
\begin{defin} (cf. \cite[page 226]{Fu94})\label{localAF}
A family $A=(A_t)_{t\ge 0}$ of functions on $\Omega$ is called an
{\it local additive functional} (local AF in short) of ${\bf M},$ if
$A$ satisfies all the requirements for an AF as stated in above (i)
and (ii), except that the additivity property (\ref{additivity}) is
 required only for $s,t\ge 0$ with
$t+s<\zeta(\omega)$.
\end{defin}
Two local AFs $A^{(1)}$, $A^{(2)}$ are said to be equivalent if for
 ${\cal E}{\textrm{-}q.e.}\ x\in E$, it holds that
$$
P_x(A^{(1)}_t=A^{(2)}_t;t<\zeta)=P_x(t<\zeta), ~~\forall~t\geq0.
$$
Define
\begin{eqnarray*}
\dot{\mathcal{M}}_{loc}
 &:=&\{M\,|\, M\ \mbox{is a local AF of}\ {\bf M},\ \exists\ \{V_n\},\{E_n\}\in\Theta\ {\rm and}\ \{M^n\,|\,M^n\in\dot{\mathcal{M}}^{V_n}\}\\
 &&\ \ \ \ \ \ \ \ \ \mbox{such that}\ E_n\subset V_n,\ M_{t\wedge\tau_{E_n}}=M^{n}_{t\wedge\tau_{E_n}},\ t\ge0,\ n\in\mathbb{N} \}
 \end{eqnarray*}
 and
 \begin{eqnarray*}
{\mathcal{N}}_{c,loc}&:=&\{N\,|\, N\ \mbox{is a local AF of}\ {\bf M},\ \exists\ \{V_n\},\{E_n\}\in\Theta\ {\rm and}\ \{N^n\,|\,N^n\in\mathcal{N}^{V_n}_{c}\}\\
 &&\ \ \ \ \ \ \ \ \ \mbox{such that}\ E_n\subset V_n,\ N_{t\wedge\tau_{E_n}}=N^{n}_{t\wedge\tau_{E_n}},\ t\ge0,\ n\in\mathbb{N} \}.
 \end{eqnarray*}
We use ${\mathcal{M}}^{[\![0, \zeta[\![}_{loc}$ to denote the family of all local martingales on ${[\![0, \zeta[\![}$ (cf. \cite[\S8.3]{HWY}).

We put the following assumption:
\begin{assum}\label{assum1}
There exists $\{{V_n}\}\in \Theta$ such that, for each $n\in\mathbb{N}$, there exists a Dirichlet form $(\eta^{(n)},D(\eta^{(n)}))$ on $L^2(V_n;m)$ and a constant $C_n>1$
such that $D(\eta^{(n)})= D(\mathcal{E})_{{V_n}}$ and for any $u\in
D(\mathcal{E})_{{V_n}}$,
\begin{eqnarray*}
\frac{1}{C_n}\eta^{(n)}_1(u,u)\leq\mathcal{E}_1(u,u)\leq C_n\eta^{(n)}_1(u,u).
\end{eqnarray*}
\end{assum}

Now we can state the main result of this section.
\begin{thm}\label{thm3.2} Suppose that
$(\mathcal{E},D(\mathcal{E}))$ is a quasi-regular local
semi-Dirichlet form on $L^{2}(E;m)$ satisfying Assumption \ref{assum1}.
Then, for any $u\in {D(\mathcal{E})}_{loc}$, there exist
$M^{[u]}\in\dot{\mathcal{M}}_{loc}$ and $N^{[u]}\in {\mathcal{N}}_{c,loc}$
such that
\begin{eqnarray}\label{new3}
\tilde{u}(X_{t})-\tilde{u}(X_{0})=M^{[u]}_{t}+N^{[u]}_{t},\ \ t\ge0,\
\ P_{x}{\textrm{-}a.s.}\ \ {\rm for}\ {\cal
E}{\textrm{-}q.e.}\ x\in E.
\end{eqnarray}
Moreover, $M^{[u]}\in{\mathcal{M}}^{[\![0, \zeta[\![}_{loc}$.

Decomposition (\ref{new3}) is unique up to the equivalence of local
AFs.
\end{thm}
Before proving Theorem \ref{thm3.2}, we present some lemmas.

We fix a $\{{V_n}\}\in\Theta$ satisfying Assumption \ref{assum1}.
Without loss of generality, we assume that $\widetilde{\hat{h}}$ is
bounded on each ${V_n}$, otherwise we may replace ${V_n}$ by
${V_n}\cap\{\widetilde{\hat{h}}< n \}$. To simplify notations,
we write
$$\bar{h}_{n}:=\bar{h}^{V_n}, ~~ \mbox{and} ~~D(\mathcal{E})_{V_n,b}:
={\cal B}_b(E) \cap D(\mathcal{E})_{V_n}.
$$
By the definition (\ref{ebn}), employing some potential theory
developed in the Appendix of this paper (cf. Lemma \ref{lnm}, Theorem
\ref{thm2.25} and Theorem \ref{thm34}), following the argument of
\cite[Theorem 5.2.1]{Fu94}, we can prove the following lemma.
\begin{lem}\label{lem3.4}
$\dot{\mathcal{M}}^{V_n}$ is a real Hilbert space with inner
product $e^{V_n}$. Moreover, if
$\{M_{l}\}\subset\dot{\mathcal{M}}^{V_n}$ is $e^{V_n}$-Cauchy,
then there exist a unique $M\in\dot{\mathcal{M}}^{V_n}$ and a
subsequence $\{l_k\}$ such that
$\lim_{k\rightarrow\infty}e^{V_n}(M_{l_k}-M)=0$ and for ${\cal
E}{\textrm{-}q.e.}\ x\in V_n$,
 $$P_{x}(\lim_{k\rightarrow\infty}M_{l_k}(t)=M(t)\ \mbox{uniformly on each
compact interval of}\ \  [0,\infty))=1.$$
 \end{lem}

Next we give Fukushima's decomposition for the part process
$X^{{V_n}}$.
\begin{lem}\label{thm3.6} Let $u\in D(\mathcal{E})_{V_n,b}$. Then
there exist unique $M^{n,[u]}\in\dot\mathcal{M}^{V_n}$ and
$N^{n,[u]}\in \mathcal{N}^{V_n}_{c}$ such that ${\rm for}\ {\cal
E}{\textrm{-}q.e.}\ x\in V_n$,
\begin{eqnarray} \label{e3.18}
\tilde{u}(X^{{V_n}}_{t})-\tilde{u}(X^{{V_n}}_{0})=M^{n,[u]}_{t}+N^{n,[u]}_{t},\
\ t\ge 0,\ \ P_{x}{\textrm{-}a.s.}
\end{eqnarray}
\end{lem}
\begin{proof} Note that if an AF $A\in\dot\mathcal{M}^{V_n}$ with $e^{V_n}(A)=0$
then $\mu^{(n)}_{<A>}(\widetilde{\bar{h}_{n}})=2e^{V_n}(A)=0$ by
Theorem \ref{thm2.25} in the Appendix  and (\ref{ebn}). Here $\mu^{(n)}_{<A>}$ denotes the Revuz measure of $A$ w.r.t. $X^{V_n}$. Hence $<A>=0$ since
$\widetilde{\bar{h}_{n}}>0$ ${\cal E}{\textrm{-}q.e.}$ on $V_n$.
Therefore $\dot\mathcal{M}^{V_n}\cap \mathcal{N}^{V_n}_{c}=\{0\}$
and the proof of the uniqueness of decomposition (\ref{e3.18}) is
complete.

To obtain the existence of decomposition (\ref{e3.18}), we start
with the special case that $u={R}_1^{V_n}f$ for some bounded Borel
function $f\in L^2(V_n;m)$, where $(R^{V_n}_t)_{t\ge 0}$ is the
resolvent of $X^{{V_n}}$. Set
\begin{eqnarray}\label{s1}
\left\{
\begin{array}{l}
N^{n,[u]}_{t}=\int^{t}_{0}(u(X^{V_n}_{s})-f(X^{V_n}_{s}))ds,\\
M^{n,[u]}_{t}=u(X^{V_n}_{t})-u(X^{V_n}_{0})-N^{n,[u]}_{t},\ \ \ \
t\ge 0.
\end{array}
\right.
\end{eqnarray}
Then $N^{n,[u]}\in \mathcal{N}^{V_n}_{c}$ and
$M^{n,[u]}\in\dot\mathcal{M}^{V_n}$. In fact,
\begin{eqnarray}\label{s2}
e^{V_n}(N^{n,[u]})
&=&\lim_{t\downarrow0}\frac{1}{2t}E_{\bar{h}_n\cdot m}[(\int^{t}_{0}(u-f)(X^{V_n}_{s})ds)^{2}]\nonumber\\
&\leq&\lim_{t\downarrow0}\frac{1}{2}E_{\bar{h}_n\cdot m}[\int^{t}_{0}(u-f)^{2}(X^{V_n}_{s})ds]\nonumber\\
&=&\lim_{t\downarrow0}\frac{1}{2}[\int^{t}_{0}\int_{V_n}\bar{h}_nT^{V_n}_{s}(u-f)^{2}dmds]\nonumber\\
&=&\lim_{t\downarrow0}\frac{1}{2}[\int^{t}_{0}\int_{V_n}(u-f)^{2}\hat{T}^{V_n}_{s}\bar{h}_ndmds]\nonumber\\
&\leq&\| u-f\|_{\infty}\lim_{t\downarrow0}\frac{1}{2}[\int^{t}_{0}\int_{V_n}|u-f|\hat{T}^{V_n}_{s}\bar{h}_ndmds]\nonumber\\
&\leq&\|
u-f\|_{\infty}\lim_{t\downarrow0}\frac{1}{2}[\int^{t}_{0}(\int_{V_n}
(u-f)^{2}dm)^{1/2}
   (\int_{V_n}(\hat{T}^{V_n}_{s}\bar{h}_n)^{2}dm)^{1/2}ds]\nonumber\\
&\leq&\| u-f\|_{\infty}(\int_{V_n}
(u-f)^{2}dm)^{1/2}(\int_{V_n}\bar{h}_n^{2}dm)^{1/2}\lim_{t\downarrow0}\frac{t}{2}\nonumber\\
&=&0.
\end{eqnarray}
By Assumption \ref{assum1}, $u^{2}\in D(\mathcal{E})_{{V_n},b}$
and $u\bar{h}_{n}\in D(\mathcal{E})_{{V_n},b}$. Then, by
(\ref{s1}), (\ref{s2}), \cite[Theorem 3.4]{f} and Assumption \ref{assum1}, we
get
  \begin{eqnarray}\label{e3.21}
 e^{V_n}\hskip -0.5cm& &\hskip -0.5cm(M^{n,[u]})\nonumber\\
 &=&\lim_{t\downarrow0}\frac{1}{2t}E_{\bar{h}_n\cdot m}[(u(X^{V_n}_{t})-u(X^{V_n}_{0}))^{2}]\nonumber\\
 \nonumber&=&\lim_{t\downarrow0}\{\frac{1}{t}(u\bar{h}_n,u-T^{V_n}_{t}
 u)-\frac{1}{2t}(\bar{h}_n,u^{2}-T^{V_n}_{t}u^{2})\}\\
   \nonumber&=&\mathcal{E}^{{V_n}}(u,u\bar{h}_{n})
      -\frac{1}{2}\mathcal{E}^{{V_n}}(u^2,\bar{h}_{n})\\
     \nonumber &\leq&\mathcal{E}^{{V_n}}_{1}(u,u\bar{h}_{n})\\
 \nonumber &\leq&
     K\mathcal{E}^{{V_n}}_{1}(u,u)^{1/2}
     \mathcal{E}_{1}^{{V_n}}(u\bar{h}_{n},u\bar{h}_{n})^{1/2}\\
 \nonumber  &\leq&
     KC_n^{1/2}\mathcal{E}_{1}^{{V_n}}(u,u)^{1/2}\eta_{1}^{{(n)}}(u\bar{h}_{n},u\bar{h}_{n})^{1/2}\\
  \nonumber&\leq&
     KC_n^{1/2}\mathcal{E}_{1}^{{V_n}}(u,u)^{1/2}(
     \|u\|_{\infty}\eta^{{(n)}}_{1}(\bar{h}_{n},\bar{h}_{n})^{1/2}
     +\|\bar{h}_{n}\|_{\infty}\eta_{1}^{{(n)}}(u,u)^{1/2})\\
 &\leq&
     KC_n\mathcal{E}_{1}^{{V_n}}(u,u)^{1/2}(
     \|u\|_{\infty}\mathcal{E}_{1}^{{V_n}}(\bar{h}_{n},\bar{h}_{n})^{1/2}
     +\|\bar{h}_{n}\|_{\infty}\mathcal{E}_{1}^{{V_n}}(u,u)^{1/2}),\ \ \ \ \ \
     \
     \end{eqnarray}
where $K$ is the continuity constant of $({\cal E},D({\cal
E}))$ (cf. (\ref{K}) in the Appendix).

Next, take any bounded Borel function $u\in D(\mathcal{E})_{V_n}$.
Define
\begin{eqnarray*}
u_{l}=lR^{V_n}_{l+1}u=R^{V_n}_{1}g_{l},\ \ \ \
g_{l}=l(u-lR^{V_n}_{l+1}u).
\end{eqnarray*}
By the uniqueness of decomposition (\ref{e3.18}) for $u_{l}$'s, we
have $M^{n, [u_l]}-M^{n,[u_{k}]}=M^{n,[u_l-u_k]}$. Then, by
(\ref{e3.21}), we get
\begin{eqnarray*}
& &e^{V_n}(M^{n, [u_l]}-M^{n,[u_{k}]})\\
&&\ \ \ \ =e^{V_n}(M^{n, [u_l-u_k]})\\
&&\ \ \ \ \leq KC_n\mathcal{E}_{1}^{{V_n}}(u_l-u_k,u_l-u_k)^{1/2}(
     \|u_l-u_k\|_{\infty}\mathcal{E}_{1}^{{V_n}}(\bar{h}_{n},\bar{h}_{n})^{1/2}\\
&&\ \ \ \ \ \ \ \ \ \ \ \ \ \ \ \
     +\|\bar{h}_{n}\|_{\infty}\mathcal{E}_{1}^{{V_n}}(u_l-u_k,u_l-u_k)^{1/2}).
\end{eqnarray*}
Since $u_{l}\in D(\mathcal{E})_{V_n}$, bounded by
$\|u\|_{\infty}$, and $\mathcal{E}^{V_n}_{1}$-convergent to $u$,
we conclude that $\{M^{n, [u_l]}\}$ is an $e^{V_n}$-Cauchy
sequence in the space
$\dot\mathcal{M}^{V_n}$. Define
$$
M^{n,[u]}=\lim_{l\rightarrow\infty}M^{n, [u_l]}\ {\rm in}\
(\dot\mathcal{M}^{V_n},e^{V_n}),\ \ \ \
N^{n,[u]}=\tilde{u}(X^{V_n}_t)-\tilde{u}(X^{V_n}_0)-M^{n,[u]}.
$$
Then $M^{n,[u]}\in \dot\mathcal{M}^{V_n}$ by Lemma \ref{lem3.4}.

It only remains to show that $N^{n,[u]}\in \mathcal{N}^{V_n}_{c}$.
By Lemma \ref{l343} in the Appendix  and Lemma \ref{lem3.4}, there exists a subsequence
$\{l_{k}\}$ such that ${\rm for}\ {\cal E}{\textrm{-}q.e.}\ x\in V_n$,
$$
P_{x}(N^{n,[u_{l_{k}}]}\ \mbox{converges to} \ N^{n,[u]}\
\mbox{uniformly  on each
compact interval of}\ \  [0,\infty))=1.
$$
From this and (\ref{s1}), we know that $N^{n,[u]}$ is a CAF. On
the other hand, by
 \begin{eqnarray*}
 N^{n,[u]}_{t}=A^{n,[u-u_l]}_{t}-(M^{n,[u]}_{t}-M^{n,[u_{l}]}_{t})+N^{n,[u_{l}]}_{t},
 \end{eqnarray*}
we get
$$
 e^{V_n}(N^{n,[u]})
  \leq3e^{V_n}(A^{n,[u-u_l]})+3e^{V_n}(M^{n,[u]}-M^{n,[u_{l}]}),
 $$
 which can be made arbitrarily small with large $l$ by
 (\ref{e3.21}). Therefore $e^{V_n}(N^{n,[u]})=0$ and $N^{n,[u]}\in \mathcal{N}^{V_n}_{c}$.
\end{proof}

We now fix a $u\in {D(\mathcal{E})}_{loc}$. Then there exist $\{V^1_n\}\in
\Theta$ and $\{u_n\}\subset D(\mathcal{E})$ such that $u=u_n$
$m{\textrm{-}a.e.}$ on $V^1_n$. By \cite[Proposition 3.6]{MR95},
we may assume without loss of generality that  each $u_n$ is
${\cal E}$-quasi-continuous. By \cite[Proposition 2.16]{MR95},
there exists an $\mathcal{E}$-nest $\{F_{n}^2\}$ of compact subsets of $E$ such that
$\{u_n\}\subset C\{F_{n}^2\}$. Denote
by $V^{2}_{n}$ the finely interior of $F^2_n$ for
$n\in\mathbb{N}$. Then $\{V^{2}_{n}\}\in\Theta$. Define
${V'_n}=V^1_n\cap V^2_{n}$. Then $\{{V'_n}\}\in\Theta$ and each
$u_n$ is bounded on ${V'_n}$. To simplify notation, we still use $V_n$ to denote $V_n\cap V'_n$ for $n\in\mathbb{N}$.

For $n\in\mathbb{N}$, we define
$E_{n}=\{x\in E\,|\,{\widetilde{h_n}}(x)>{1\over n}\}$, where
$h_n:=G_{1}^{{V_n}}\phi$. Then $\{E_{n}\}\in\Theta$ satisfying
$\overline{E}_n^{\mathcal{E}}\subset E_{n+1}\ {\cal E}{\textrm{-}q.e.}$
and $E_{n}\subset {V_n}\ {\cal E}{\textrm{-}q.e.}$ for each $n\in
\mathbb{N}$ (cf. \cite[Lemma 3.8]{kuwae}). Here $\overline{E}_n^{\mathcal{E}}$ denotes the ${\cal E}$-quasi-closure of $E_n$. Define
$f_{n}=n\widetilde{h_n}\wedge1$. Then $f_{n}=1$ on $E_{n}$ and
$f_{n}=0$ on $V^c_{n}$. Since $f_{n}$ is a 1-excessive function of
$(\mathcal{E}^{V_n},D(\mathcal{E})_{V_n})$ and $f_{n}\leq
n\widetilde{h_n}\in D(\mathcal{E})_{V_n}$, hence $f_{n}\in
D(\mathcal{E})_{V_n}$ by \cite[Remark 3.4(ii)]{MR1995}. Denote by
$Q_n$ the bound of $|u_n|$ on $V_n$. Then  $u_{n}f_{n}=((-Q_n)\vee
u_n\wedge Q_n)f_n\in D(\eta)_{V_n,b}=D(\mathcal{E})_{V_n,b}$.

For $n\in \mathbb{N}$, we denote by $\{\mathcal{F}^n_{t}\}$ the
minimum completed admissible filtration of $X^{{V_n}}$. For $n<l$,  $\mathcal{F}^n_{t}\subset \mathcal{F}^l_{t}\subset
\mathcal{F}_{t}$. Since $E_n\subset V_n$, $\tau_{E_n}$ is an
$\{\mathcal{F}^{n}_{t}\}$-stopping time.

\begin{lem}\label{lem3.7} For $n<l$, we have
 $M^{n,[u_nf_n]}_{t\wedge\tau_{E_n}}=M^{l,[u_lf_l]}_{t\wedge\tau_{E_n}}$
 and $N^{n,[u_nf_n]}_{t\wedge\tau_{E_n}}=N^{l,[u_lf_l]}_{t\wedge\tau_{E_n}}$, $t\ge0$, $P_{x}{\textrm{-}a.s.}\ \ {\rm for}\ {\cal
E}{\textrm{-}q.e.}\ x\in  V_n$.
\end{lem}
\begin{proof} Let $n<l$. Since
$M^{n,[u_nf_n]}\in\dot\mathcal{M}^{V_n}$, $M^{n,[u_nf_n]}$ is an
$\{\mathcal{F}^{n}_{t}\}$-martingale by the Markov property. Since
$\tau_{E_n}$ is an $\{\mathcal{F}^{n}_{t}\}$-stopping time,
$\{M^{n,[u_nf_n]}_{_{t\wedge\tau_{E_{n}}}}\}$ is an
$\{\mathcal{F}^{n}_{t\wedge\tau_{E_{n}}}\}$-martingale. Denote
$\Upsilon^n_t=\sigma\{X^{{V_n}}_{s\wedge\tau_{E_n}}\,|\,0\leq
s\leq t\}$. Then $\{M^{n,[u_nf_n]}_{_{t\wedge\tau_{E_{n}}}}\}$ is
a $\{\Upsilon^n_{t}\}$-martingale. Denote
$\Upsilon^{n,l}_t=\sigma\{X^{{V_l}}_{s\wedge\tau_{E_n}}\,|\,0\leq
s\leq t\}$. Similarly, we can show that
$\{M^{l,[u_nf_n]}_{_{t\wedge\tau_{E_{n}}}}\}$ is a
$\{\Upsilon^{n,l}_{t}\}$-martingale. By the assumption that ${\bf M}$ is a diffusion, the fact that $f_n$ is quasi-continuous and $f_n=1$ on $E_n$, we get $f_n(X_{s\wedge\tau_{E_n}})=1$ if $0<s\wedge\tau_{E_n}<\zeta$. Hence $X_{s\wedge\tau_{E_n}}\in V_n$, if $0<s\wedge\tau_{E_n}<\zeta$, since $f_n=0$ on $V_n^c$. Therefore
\begin{equation}\label{eqadd}
X^{V_l}_{s\wedge\tau_{E_n}}
 =X_{s\wedge\tau_{E_n}} =X_{s\wedge\tau_{E_n}}^{{V_n}},\ \ P_{x}{\textrm{-}a.s.}\ \ {\rm for}\ {\cal
E}{\textrm{-}q.e.}\ x\in V_n,
\end{equation}
which implies that $\{M^{l,[u_nf_n]}_{_{t\wedge\tau_{E_{n}}}}\}$ is a
$\{\Upsilon^n_t\}$-martingale.

Let $N\in \mathcal{N}^{V_j}_c$ for some $j\in
\mathbb{N}$. Then, for any $T>0$,
\begin{eqnarray*}
\sum_{k=1}^{[rT]}E_{\bar{h}_j\cdot
m}[(N_{\frac{k+1}{r}}-N_{\frac{k}{r}})^{2}]
&\leq&\sum_{k=1}^{[rT]}e^{T}(E_{\cdot}(N_{\frac{1}{r}}^{2}),e^{-\frac{k}{r}}\hat{T}^{V_j}_{\frac{k}{r}}\bar{h}_j)\\
&\leq&\sum_{k=1}^{[rT]}e^{T}(E_{\cdot}(N_{\frac{1}{r}}^{2}),\bar{h}_j)\\
&\leq&rTe^{T}E_{\bar{h}_j\cdot
m}(N_{\frac{1}{r}}^{2})\rightarrow0\ \ \ \mbox{as} \ \
r\rightarrow\infty.
\end{eqnarray*}
Hence
$$\sum_{k=1}^{[rT]}(N_{\frac{k+1}{r}}-N_{\frac{k}{r}})^{2}\rightarrow0,
\ \ r\rightarrow\infty,\ \ {\rm in}\ \ P_{m},$$ which implies that the
quadratic variation process of $N$ w.r.t. $P_m$ is 0.

By \cite[Proposition 3.3]{kuwae},
 $(\widehat{\hat{G}_{1}\phi})^{1}_{V^{c}_{n}}=\hat{G}_{1}\phi-\hat{G}_{1}^{{V_n}}\phi$.
Since $V^{c}_{n}\supset V^{c}_l$,
$(\widehat{\hat{G}_{1}\phi})^{1}_{V^{c}_{n}}\geq
(\widehat{\hat{G}_{1}\phi})^{1}_{V^{c}_l}$. Then
$\hat{G}_{1}^{{V_n}}\phi\leq \hat{G}_{1}^{V_l}\phi$ and thus
\begin{equation}\label{new2}\bar{h}_{n}\le\bar{h}_{l}.
\end{equation}
Therefore
\begin{equation}\label{new1} e^{V_n}(A)\le e^{V_l}(A)
\end{equation}
for any AF $A=(A_t)_{t\ge 0}$ of $X^{{V_n}}$.

Note that $N^{l,[u_nf_n]}_{t\wedge\tau_{E_{n}}}=(\widetilde{u_nf_n})(X^{V_l}_{t\wedge\tau_{E_{n}}})
-(\widetilde{u_nf_n})(X^{V_l}_{0})-M^{l,[u_nf_n]}_{t\wedge\tau_{E_{n}}}\in\Upsilon^{n,l}_t=
\Upsilon^{n}_t\subset\mathcal{F}^n_{t\wedge\tau_{E_{n}}}$. By the analog of \cite[Lemma 5.5.2]{Fu94} in the semi-Dirichlet forms setting,
$\{N^{l,[u_nf_n]}_{t\wedge\tau_{E_{n}}}\}$ is a CAF of
$X^{{V_n}}$. By (\ref{new1}), $e^{V_n}(N^{l,[u_nf_n]}_{t\wedge\tau_{E_{n}}})\leq
e^{V_l}(N^{l,[u_nf_n]}_{t\wedge\tau_{E_{n}}})=0$. Hence $(N^{l,[u_nf_n]}_{t\wedge\tau_{E_{n}}})_{t\ge 0}
\in \mathcal{N}^{V_n}_c$, which implies that
the quadratic variation process of $\{N^{l,[u_nf_n]}_{t\wedge\tau_{E_{n}}}\}$ w.r.t. $P_m$
is 0. Since ${\rm for}\ {\cal E}{\textrm{-}q.e.}\ x\in V_n$, by (\ref{eqadd}),
\begin{eqnarray*}
M^{n,[u_nf_n]}_{t\wedge\tau_{E_{n}}}+N^{n,[u_nf_n]}_{t\wedge\tau_{E_{n}}}
&=&\widetilde{u_nf_n}(X^{{V_n}}_{t\wedge\tau_{E_n}})-\widetilde{u_nf_n}(X^{{V_n}}_0)\\
&=&\widetilde{u_nf_n}(X^{V_l}_{t\wedge\tau_{E_n}})-\widetilde{u_nf_n}(X^{V_l}_0)\\
&=&M^{l,[u_nf_{n}]}_{t\wedge\tau_{E_{n}}}+N^{l,[u_nf_{n}]}_{t\wedge\tau_{E_{n}}},
\ \ P_x-a.s.,
\end{eqnarray*}
and both $\{M^{n,[u_nf_n]}_{_{t\wedge\tau_{E_{n}}}}\}$ and $\{M^{l,[u_nf_n]}_{_{t\wedge\tau_{E_{n}}}}\}$ are
$\{\Upsilon^n_t\}$-martingale, hence $M^{n,[u_nf_n]}_{t\wedge\tau_{E_{n}}}=M^{l,[u_nf_n]}_{t\wedge\tau_{E_{n}}}$
 and $N^{n,[u_nf_n]}_{t\wedge\tau_{E_{n}}}=N^{l,[u_nf_n]}_{t\wedge\tau_{E_{n}}}$, $P_{x}{\textrm{-}a.s.}\ \ {\rm for}\ m{\textrm{-}a.e.}\ x\in V_n$. This implies that
$E_m(<M^{n,[u_nf_n]}_{\cdot\wedge\tau_{E_{n}}}-M^{l,[u_nf_n]}_{\cdot\wedge\tau_{E_{n}}}>_t)=0$, $\forall t\ge 0$. Then, by Theorem \ref{thm2.25}(i) in the Appendix, $M^{n,[u_nf_n]}_{t\wedge\tau_{E_{n}}}=M^{l,[u_nf_n]}_{t\wedge\tau_{E_{n}}}$, $\forall t\ge 0$,
 $P_{x}{\textrm{-}a.s.}\ \ {\rm for}\ {\cal
E}{\textrm{-}q.e.}\ x\in V_n$. Hence $N^{n,[u_nf_n]}_{t\wedge\tau_{E_{n}}}=N^{l,[u_nf_n]}_{t\wedge\tau_{E_{n}}}$, $\forall t\ge 0$, $P_{x}{\textrm{-}a.s.}\ \ {\rm for}\ {\cal
E}{\textrm{-}q.e.}\ x\in V_n$.

Since $u_nf_n=u_{l}f_{l}=u$ on $E_n$, similar to \cite[Lemma 2.4]{kuwae2010}, we can show that
$M^{l,[u_{n}f_{n}]}_{t}=M^{l,[u_lf_l]}_{t}$ when $t<\tau_{E_{n}}$, $P_{x}{\textrm{-}a.s.}\ \ {\rm for}\ {\cal
E}{\textrm{-}q.e.}\ x\in V_l$. If $\tau_{E_{n}}=\zeta$, then by the fact $u_{n}f_{n}(X^{V_l}_{\zeta})=u_{l}f_{l}(X^{V_l}_{\zeta})=0$ and the continuity of $N^{l,[u_{n}f_{n}]}_{t}$ and $N^{l,[u_lf_l]}_{t}$, one finds that $M^{l,[u_{n}f_{n}]}_{t\wedge\tau_{E_{n}}}=M^{l,[u_lf_l]}_{t\wedge\tau_{E_{n}}}$. By the quasi-continuity of $u_nf_n$, $u_{l}f_{l}$ and the assumption that ${\bf M}$ is a diffusion, one finds that $M^{l,[u_{n}f_{n}]}$ and $M^{l,[u_lf_l]}$ are continuous on $[0,\zeta)$, $P_{x}{\textrm{-}a.s.}\ \ {\rm for}\ {\cal
E}{\textrm{-}q.e.}\ x\in V_l$. Hence, if $\tau_{E_{n}}<\zeta$ we have $M^{l,[u_{n}f_{n}]}_{\tau_{E_{n}}}=M^{l,[u_lf_l]}_{\tau_{E_{n}}}$. Therefore $M^{n,[u_nf_n]}_{t\wedge\tau_{E_n}}=M^{l,[u_lf_l]}_{t\wedge\tau_{E_n}}$
 and $N^{n,[u_nf_n]}_{t\wedge\tau_{E_n}}=N^{l,[u_lf_l]}_{t\wedge\tau_{E_n}}$, $t\ge0$, $P_{x}{\textrm{-}a.s.}\ \ {\rm for}\ {\cal
E}{\textrm{-}q.e.}\ x\in  V_n$.
\end{proof}

\noindent\textbf{Proof of Theorem \ref{thm3.2}}\ \
We define $M^{[u]}_{t\wedge\tau_{E_n}}:=\lim_{l\rightarrow\infty}M^{l,[u_lf_l]}_{t\wedge\tau_{E_n}}$ and $M^{[u]}_t:=0$ for $t>\zeta$ if there exists some $n$ such that $\tau_{E_n}=\zeta$ and $\zeta<\infty$; or $M^{[u]}_t:=0$ for $t\ge\zeta$, otherwise. By Lemma \ref{lem3.7}, $M^{[u]}$ is well defined. Define $M^n_t:=M^{n+1,[u_{n+1}f_{n+1}]}_{t\wedge\tau_{E_n}}$ for $t\ge 0$ and $n\in\mathbb{N}$. Then $M^{[u]}_{t\wedge\tau_{E_n}}=M^{n}_{t\wedge\tau_{E_n}}$ $P_{x}{\textrm{-}a.s.}\ \ {\rm for}\ {\cal
E}{\textrm{-}q.e.}\ x\in  V_{n+1}$ by Lemma \ref{lem3.7}. Since $\overline{E}_n^{\mathcal{E}}\subset E_{n+1}\subset V_{n+1}\ {\cal E}{\textrm{-}q.e.}$ implies that $P_x(\tau_{E_n}=0)=1$ for $x\notin V_{n+1}$, $M^{[u]}_{t\wedge\tau_{E_n}}=M^{n}_{t\wedge\tau_{E_n}}$ $P_{x}{\textrm{-}a.s.}\ \ {\rm for}\ {\cal
E}{\textrm{-}q.e.}\ x\in  E$. Similar to (\ref{new2}) and (\ref{new1}), we can show that  $e^{V_n}(M^{n})\le e^{V_{n+1}}(M^{n})$ for each $n\in\mathbb{N}$. Then $M^{n}\in \dot{\mathcal{M}}^{V_n}$ and hence $M^{[u]}\in \dot{\mathcal{M}}_{loc}$. Define $N^{[u]}_t=\tilde{u}(X_{t})-\tilde{u}(X_{0})-M^{[u]}_t$. Then, we have $N^{[u]}_{t\wedge\tau_{E_n}}=\lim_{l\rightarrow\infty}N^{l,[u_lf_l]}_{t\wedge\tau_{E_n}}$. Moreover $N^{[u]}\in {\mathcal{N}}_{c,loc}$.

Next we show that $M^{n}$ is also an $\{\mathcal{F}_{t}\}$-martingale, which implies that $M^{[u]}\in{\mathcal{M}}^{[\![0, \zeta[\![}_{loc}$. In fact, by the fact that $\tau_{E_n}$ is an $\{\mathcal{F}^{n+1}_{t}\}$-stopping time, we find that $I_{\tau_{E_n}\le s}$ is
${\cal F}^{n+1}_{s\wedge\tau_{E_n}}$-measurable for any $s\ge 0$. Let $0\le s_1<\cdots<s_k\le s<t$ and $g\in{\cal B}_b(\mathbb{R}^k)$. Then, we obtain by (\ref{eqadd}) and the fact $M^{n+1,[u_{n+1}f_{n+1}]}\in \dot{\mathcal{M}}^{V_{n+1}}$ that ${\rm for}\ {\cal
E}{\textrm{-}q.e.}\ x\in  V_{n+1}$,
\begin{eqnarray*}
& &\int_{\Omega} M^n_tg(X_{s_1},\dots,X_{s_k})dP_x\\&&\ \ \ \ =\int_{\tau_{E_n}\le s}M^n_tg(X_{s_1},\dots,X_{s_k})dP_x
+\int_{\tau_{E_n}>s}M^n_tg(X_{s_1},\dots,X_{s_k})dP_x\\
&&\ \ \ \ =\int_{\tau_{E_n}\le s}M^n_sg(X_{s_1},\dots,X_{s_k})dP_x\\
&&\ \ \ \ \ \ \ \ +\int_{\Omega} M^{n+1,[u_{n+1}f_{n+1}]}_{t\wedge\tau_{E_n}}g(X^{V_{n+1}}_{s_1\wedge\tau_{E_n}},
\dots,X^{V_{n+1}}_{s_k\wedge\tau_{E_n}})
I_{\tau_{E_n}> s}dP_x\\
&&\ \ \ \ =\int_{\tau_{E_n}\le s}M^n_sg(X_{s_1},\dots,X_{s_k})dP_x\\
&&\ \ \ \ \ \ \ \ +\int_{\Omega} M^{n+1,[u_{n+1}f_{n+1}]}_{s\wedge\tau_{E_n}}g(X^{V_{n+1}}_{s_1\wedge\tau_{E_n}},
\dots,X^{V_{n+1}}_{s_k\wedge\tau_{E_n}})
I_{\tau_{E_n}> s}dP_x\\
&&\ \ \ \ =\int_{\tau_{E_n}\le s}M^n_sg(X_{s_1},\dots,X_{s_k})dP_x
+\int_{\tau_{E_n}>s}M^n_sg(X_{s_1},\dots,X_{s_k})dP_x\\
&&\ \ \ \ =\int_{\Omega} M^n_sg(X_{s_1},\dots,X_{s_k})dP_x.
\end{eqnarray*}
Obviously, the equality holds for $x\notin  V_{n+1}$. Therefore,
$M^{n}$ is an $\{\mathcal{F}_{t}\}$-martingale.

Finally, we prove the uniqueness of decomposition (\ref{new3}).
Suppose that $M^1\in\dot\mathcal{M}_{loc}$ and
$N^1\in{\mathcal{N}}_{c,loc}$ such that
$$
\tilde{u}(X_{t})-\tilde{u}(X_{0})=M^1_{t}+N^1_{t},\ \ t\ge0,\ \ P_{x}{\textrm{-}a.s.}\ \ {\rm for}\ {\cal
E}{\textrm{-}q.e.}\ x\in E.
$$
Then, there exists $\{E_n\}\in\Theta$ such that, for each $n\in\mathbb{N}$, $\{(M^{[u]}-M^1)I_{[\![0, \tau_{E_n}]\!]}\}$ is a square integrable martingale and a zero quadratic variation process w.r.t. $P_m$. This implies that
$P_m(<(M^{[u]}-M^1)I_{[\![0, \tau_{E_n}]\!]}>_t=0, \forall t\in[0,\infty))=0$. Consequently by the analog of \cite[Lemma 5.1.10]{Fu94} in the semi-Dirichlet forms setting, $P_x(<(M^{[u]}-M^1)I_{[\![0, \tau_{E_n}]\!]}>_t=0, \forall t\in[0,\infty))=0$ ${\rm for}\ {\cal
E}{\textrm{-}q.e.}\ x\in E$. Therefore $M_t^{[u]}=M_t^1$, $0\le t\le\tau_{E_n}$, $P_{x}{\textrm{-}a.s.}\ {\rm for}\ {\cal
E}{\textrm{-}q.e.}\ x\in E$. Since $n$ is arbitrary, we obtain the uniqueness of decomposition (\ref{new3}) up to the equivalence of local AFs.\hfill\fbox

 %%%%%%%%%%%%%%%%%%%%%%%%%%%%%%%%%%%%%%%%%%%%%%%%%%%%%%%%%%%%%%%%%%%%%%%%%%%%%%%%%%%%%%%%%%%%%%%%%%%%%%%%%%%%%%%%%
\section {Transformation formula}\setcounter{equation}{0}
\label{sec:transform}In this section, we adopt the setting of
Section 2. Suppose that $(\mathcal{E},D(\mathcal{E}))$ is a
quasi-regular local semi-Dirichlet form on $L^{2}(E;m)$ satisfying
Assumption \ref{assum1}. We fix a $\{{V_n}\}\in\Theta$ satisfying
Assumption \ref{assum1} and satisfying that $\widetilde{\hat{h}}$ is
bounded on each $V_n$. Let $X^{V_n}$,
$(\mathcal{E}^{V_n},D(\mathcal{E})_{V_n})$, $\bar{h}_n$, etc. be the
same as in Section 2.
 For $u\in D(\mathcal{E})_{V_n,b}$, we denote by $\mu^{(n)}_{<u>}$ the Revuz
 measure of $<M^{n,[u]}>$ (cf. Lemma \ref{thm3.6} and Theorem \ref{thm2.25} in the Appendix). For $u,v\in D(\mathcal{E})_{V_n,b}$, we define
\begin{equation}\label{p01}
 \mu^{(n)}_{<u,v>}:=\frac{1}{2}(\mu^{(n)}_{<u+v>}-\mu^{(n)}_{<u>}-\mu^{(n)}_{<v>}).
\end{equation}

 \begin{lem}\label{lem4.1}
 Let $u,v,f\in D(\mathcal{E})_{V_n,b}$. Then
 \begin{eqnarray}\label{e4.27*}
 \int_{V_n}\tilde{f}d\mu^{(n)}_{<u,v>}=\mathcal{E}(u,vf)+\mathcal{E}(v,uf)
 -\mathcal{E}(uv,f).
 \end{eqnarray}
 \end{lem}
 \begin{proof}  By the polarization identity, (\ref{e4.27*}) holds for $u,v,f\in D(\mathcal{E})_{V_n,b}$ is equivalent to
 \begin{eqnarray}\label{au2}
 \int_{V_n}\tilde{f}d\mu^{(n)}_{<u>}=2\mathcal{E}(u,uf)-\mathcal{E}(u^2,f),\ \ \forall u,f\in D(\mathcal{E})_{V_n,b}.
 \end{eqnarray}
Below, we will prove (\ref{au2}). Without loss of generality, we assume that $f\geq0$.

For $k,l\in \mathbb{N}$, we define $f_{k}:=f\wedge (k\bar{h}_n)$ and $f_{k,l}:=l\hat{G}^{V_n}_{l+1}f_k$. By \cite[(3.9)]{MR1995}, $f_{k}\in D(\mathcal{E})_{V_n,b}$ and
\begin{equation}\label{au11}
{\cal E}_1(f_k,f_k)\le {\cal E}_1(f,f_k).
\end{equation}
By \cite[Proposition III.1.2]{MR92}, $f_{k,l}$ is $(l+1)$-co-excessive. Since $\bar{h}_n$ is 1-co-excessive,
\begin{equation}\label{au4}
0\le f_{k,l}\le k\bar{h}_n.
\end{equation} Hence $f_{k,l}\in D(\mathcal{E})_{V_n,b}$ by noting that $\bar{h}_n$ is bounded.

Note that by (\ref{au4})
 \begin{equation}\label{au1}
\lim_{t\downarrow0}{1\over t}E_{f_{k,l}\cdot m}[(N^{n,[u]}_t)^2]
 \leq k\lim_{t\downarrow0}{1\over t}E_{\bar{h}_n\cdot m}[(N^{n,[u]}_t)^2]
 =2ke^{V_n}(N^{n,[u]})=0.
 \end{equation}
Then, by Theorem \ref{thm2.25}(i) in the Appendix and (\ref{au1}), we get
  \begin{eqnarray}\label{au3}
 \int_{V_n}\widetilde{f_{k,l}}d\mu^{(n)}_{<u>}&=&\lim_{t\downarrow0}{1\over t}E_{f_{k,l}\cdot m}[<M^{n,[u]}>_t]\nonumber\\
  &=&\lim_{t\downarrow0}{1\over t}E_{f_{k,l}\cdot m}[(\tilde{u}(X^{V_n}_t)-\tilde{u}(X^{V_n}_0))^2]\nonumber\\
  &=&\lim_{t\downarrow0}{2\over t}(uf_{k,l},u-P^{V_n}_tu)
     -\lim_{t\downarrow0}{1\over t}(f_{k,l},u^2-P^{V_n}_{t}u^2)\nonumber\\
  &=&2\mathcal{E}(u,uf_{k,l}) -\mathcal{E}(u^2,f_{k,l}).
  \end{eqnarray}

By \cite[Theorem I.2.13]{MR92}, for each $k\in \mathbb{N}$, $f_{k,l}\rightarrow f_k$ in $D(\mathcal{E})_{V_n}$ as $l\rightarrow\infty$. Furthermore, by Assumption \ref{assum1}, \cite[Corollary I.4.15]{MR92} and (\ref{au4}), we can show that
$\sup_{l\geq1}\mathcal{E}(uf_{k,l},uf_{k,l})<\infty$. Thus, we obtain by \cite[Lemma I.2.12]{MR92} that $uf_{k,l}\rightarrow uf_k$
weakly in $D(\mathcal{E})_{V_n}$ as $l\rightarrow\infty$. Note that $\int_{V_n}\widetilde{\bar{h}_n}d\mu^{(n)}_{<u>}=2e^{V_n}(M^{n,[u]})<\infty$ for any $u\in D(\mathcal{E})_{V_n,b}$.
Therefore, we obtain by (\ref{au3}), (\ref{au4}) and the dominated convergence theorem that
\begin{eqnarray}\label{au5}
 \int_{V_n}\widetilde{f_k}d\mu^{(n)}_{<u>}=2\mathcal{E}(u,uf_k)-\mathcal{E}(u^2,f_k),\ \ \forall u\in D(\mathcal{E})_{V_n,b}.
 \end{eqnarray}

By (\ref{au11}) and the weak sector condition, we get
$\sup_{k\geq1}\mathcal{E}_{1}(f_k,f_k)<\infty$. Furthermore, by Assumption \ref{assum1} and \cite[Corollary I.4.15]{MR92}, we can show that
$\sup_{k\geq1}\mathcal{E}(uf_k,uf_k)<\infty$. Thus, we obtain by \cite[Lemma I.2.12]{MR92} that $f_k\rightarrow f$ and $uf_k\rightarrow uf$
weakly in $D(\mathcal{E})_{V_n}$ as $k\rightarrow\infty$. Therefore (\ref{au2}) holds by (\ref{au5}) and the monotone convergence theorem.
  \end{proof}

For $u\in {D(\mathcal{E})}_{V_n,b}$, we denote by $M^{n,[u],c}$ and $M^{n,[u],k}$ the continuous and killing parts of $M^{n,[u]}$, respectively; denote by $\mu^{n,c}_{<u>}$ and $\mu^{n,k}_{<u>}$ the Revuz
 measures of $<M^{n,[u],c}>$ and $<M^{n,[u],k}>$, respectively. Then
$M^{n,[u]}=M^{n,[u],c}+M^{n,[u],k}$ with
$$M^{n,[u],k}=-\tilde{u}(X^{V_n}_{\zeta^{(n)}-})I_{\{\zeta^{(n)}\le t\}}-(-\tilde{u}(X^{V_n}_{\zeta^{(n)}-})I_{\{\zeta^{(n)}\le t\}})^p,$$
where $\zeta^{(n)}$ denotes the life time of $X^{V_n}$ and $p$ denotes the dual predictable projection, and
\begin{equation}\label{n00}
\mu^{(n)}_{<u>}=\mu^{n,c}_{<u>}+\mu^{n,k}_{<u>}.
\end{equation}
 Let $(N^{(n)}(x,dy),H^{(n)})$ be a L\'evy system of $X^{V_n}$ and $\nu^{(n)}$ be the Revuz measure of $H^{(n)}$. Define $K^{(n)}(dx):=N^{(n)}(x,\Delta)\nu^{(n)}(dx)$. Similar to \cite[(5.3.8) and (5.3.10)]{Fu94}, we can show that
  \begin{eqnarray}\label{l00}
<M^{n,[u],k}>_t&=&(\tilde{u}^2(X^{V_n}_{\zeta^{(n)}-})I_{\zeta^{(n)}\le t})^p\nonumber\\
&=&\int_0^t\tilde{u}^2(X^{V_n}_s)N^{(n)}(X^{V_n}_s,\Delta)dH^{(n)}_s
  \end{eqnarray}
  and
    \begin{equation}\label{xc1}
\mu^{n,k}_{<u>}(dx)=\tilde{u}^2(x)K^{(n)}(dx).
\end{equation}
   For $u,v\in D(\mathcal{E})_{V_n,b}$, we define
\begin{equation}\label{n02}
 \mu^{n,c}_{<u,v>}:=\frac{1}{2}(\mu^{n,c}_{<u+v>}-\mu^{n,c}_{<u>}-\mu^{n,c}_{<v>}),\ \ \mu^{n,k}_{<u,v>}:=\frac{1}{2}(\mu^{n,k}_{<u+v>}-\mu^{n,k}_{<u>}-\mu^{n,k}_{<v>}).
\end{equation}

\begin{thm}\label{th4.5}
Let $u,v,w\in D(\mathcal{E})_{V_n,b}$. Then
\begin{eqnarray}\label{m00}
d\mu^{n,c}_{<uv,w>}=\tilde{u}d\mu^{n,c}_{<v,w>}+\tilde{v}d\mu^{n,c}_{<u,w>}.
\end{eqnarray}
\end{thm}
\begin{proof} By quasi-homeomorphism and the polarization identity, (\ref{m00}) holds for $u,v,w\in D(\mathcal{E})_{V_n,b}$ is equivalent to
\begin{eqnarray}\label{e4.28*}
\int_{V_n} \tilde{f}d\mu^{n,c}_{<u^2,w>}=2\int_{V_n} \tilde{f}\tilde{u}d\mu^{n,c}_{<u,w>},\ \ \forall f,u,w\in D(\mathcal{E})_{V_n,b}.
\end{eqnarray}
By (\ref{p01}) and (\ref{n00})-(\ref{n02}), we find that (\ref{e4.28*}) is equivalent to
\begin{eqnarray}\label{p00}
\int_{V_n} \tilde{f}d\mu^{(n)}_{<u^2,w>}+\int_{V_n} \tilde{f}\tilde{u}^2\tilde{w}dK^{(n)}=2\int_{V_n} \tilde{f}\tilde{u}d\mu^{(n)}_{<u,w>},\ \ \forall f,u,w\in D(\mathcal{E})_{V_n,b}.\ \
\end{eqnarray}

For $k\in \mathbb{N}$, we define $v_k:=kR^{V_n}_{k+1}u$. Then $v_k\rightarrow u$ in $D(\mathcal{E})_{V_n}$ as $k\rightarrow\infty$.
By Assumption \ref{assum1} and \cite[Corollary I.4.15]{MR92}, we can show that
$\sup_{k\geq1}\mathcal{E}(v_kw,v_kw)<\infty$. Then, by \cite[Lemma I.2.12]{MR92}, there exists a subsequence $\{(v_{k_l})\}_{l\in\mathbb{N}}$ of $\{v_k\}_{k\in\mathbb{N}}$ such that $u_kw\rightarrow uw$ in $D(\mathcal{E})_{V_n}$ as $k\rightarrow\infty$, where $u_k:=\frac{1}{k}\sum_{l=1}^kv_{k_l}$. Note that $u_k\rightarrow u$ in $D(\mathcal{E})_{V_n}$ as $k\rightarrow\infty$  and $\|u_k\|_{\infty}\le \|u\|_{\infty}$ for $k\in \mathbb{N}$. Moreover, $\|L^{V_n}u_k\|_{\infty}<\infty$ for $k\in \mathbb{N}$, where $L^{V_n}$ is the generator of $X^{V_n}$.

By Assumption \ref{assum1} and \cite[Corollary I.4.15]{MR92}, we can show that
$\sup_{k\geq1}[\mathcal{E}(u_kfw,u_kfw)+\mathcal{E}(u^2_kf,u^2_kf)+\mathcal{E}(u_kf,u_kf)]<\infty$. Then, we obtain by \cite[Lemma I.2.12]{MR92} that $u_kfw\rightarrow ufw$, $u^2_kf\rightarrow u^2f$ and $u_kf\rightarrow uf$
weakly in $D(\mathcal{E})_{V_n}$ as $k\rightarrow\infty$. Hence by  (\ref{e4.27*}) and  the fact $\sup_{k\ge 1}[\mathcal{E}(u_kfw,u_kfw)+\mathcal{E}(u_kf,u_kf)]<\infty$ we get
\begin{eqnarray}\label{p10}
\int_{V_n} \tilde{f}\tilde{u}d\mu^{(n)}_{<u,w>}&=&\mathcal{E}(u,ufw)+\mathcal{E}(w,u^2f)
 -\mathcal{E}(uw,uf)\nonumber\\
 &=&\lim_{k\rightarrow\infty}[\mathcal{E}(u,u_kfw)+\mathcal{E}(w,u_k^2f)
 -\mathcal{E}(uw,u_kf)]\nonumber\\
 &=&\lim_{k\rightarrow\infty}[\mathcal{E}(u_k,u_kfw)+\mathcal{E}(w,u_k^2f)
 -\mathcal{E}(u_kw,u_kf)]\nonumber\\
 &=&\lim_{k\rightarrow\infty}\int_{V_n} \tilde{f}\widetilde{u_k}d\mu^{(n)}_{<u_k,w>}.
\end{eqnarray}
By Assumption \ref{assum1} and \cite[Corollary I.4.15]{MR92}, we can show that
$\sup_{k\geq1}[\mathcal{E}(u^2_k,u^2_k)+\mathcal{E}(u^2_kf,u^2_kf)+\mathcal{E}(u^2_kw,u^2_kw)]<\infty$. Then, we obtain by \cite[Lemma I.2.12]{MR92} that $u^2_k\rightarrow u^2$, $u^2_kf\rightarrow u^2f$ and $u^2_kw\rightarrow u^2w$
weakly in $D(\mathcal{E})_{V_n}$ as $k\rightarrow\infty$. Hence by  (\ref{e4.27*}) we get
\begin{eqnarray}\label{p11}
\int_{V_n} \tilde{f}d\mu^{(n)}_{<u^2,w>}&=&\mathcal{E}(u^2,fw)+\mathcal{E}(w,u^2f)
 -\mathcal{E}(u^2w,f)\nonumber\\
 &=&\lim_{k\rightarrow\infty}[\mathcal{E}(u_k^2,fw)+\mathcal{E}(w,u_k^2f)
 -\mathcal{E}(u_k^2w,f)]\nonumber\\
&=&\lim_{k\rightarrow\infty}\int_{V_n} \tilde{f}d\mu^{(n)}_{<u^2_k,w>}.
\end{eqnarray}
By (\ref{p10}), (\ref{p11}) and the dominated convergence theorem, to prove (\ref{p00}), we may assume without loss of generality that $u$ is equal to some $u_k$. Moreover, we assume without loss of generality that $f\geq0$.

For $k,l\in \mathbb{N}$, we define $f_{k}:=f\wedge (k\bar{h}_n)$ and $f_{k,l}:=l\hat{G}^{V_n}_{l+1}f_k$. By \cite[(3.9)]{MR1995}, $f_{k}\in D(\mathcal{E})_{V_n,b}$; by \cite[Proposition III.1.2]{MR92}, $f_{k,l}$ is $(l+1)$-co-excessive. Since $\bar{h}_n$ is 1-co-excessive,
$$
0\le f_{k,l}\le k\bar{h}_n.
$$ Hence $f_{k,l}\in D(\mathcal{E})_{V_n,b}$ by noting that $\bar{h}_n$ is bounded. By the dominated convergence theorem, to prove that (\ref{p00}) holds for any $f\in D(\mathcal{E})_{V_n,b}$, it suffices to prove that (\ref{p00}) holds for any $f_{k,l}$.

\emph{Below, we will prove (\ref{p00}) for $u=u_k$ and $f=f_{k,l}$.}

Note that for any $g\in D(\mathcal{E})_{V_n,b}$,
 \begin{equation}\label{aum1}
\lim_{t\downarrow0}{1\over t}E_{f_{k,l}\cdot m}[(N^{n,[g]}_t)^2]
 \leq k\lim_{t\downarrow0}{1\over t}E_{\bar{h}_n\cdot m}[(N^{n,[g]}_t)^2]
 =2ke^{V_n}(N^{n,[g]})=0.
 \end{equation}
By Theorem \ref{thm2.25}(i) in the Appendix and (\ref{aum1}), we get
\begin{eqnarray}\label{q00}
\int_{V_n} \widetilde{f_{k,l}}d\mu^{(n)}_{<u_k^2,w>}&=&\lim_{t\downarrow0}{1\over t}E_{f_{k,l}\cdot m}[<M^{n,[u_k^2]},M^{n,[w]}>_t]\nonumber\\
&=&\lim_{t\downarrow0}\frac{1}{t}E_{f_{k,l}\cdot m}[(\widetilde{u_k}^2(X^{V_n}_t)-\widetilde{u_k}^2(X^{V_n}_0))(\tilde{w}(X^{V_n}_t)-\tilde{w}(X^{V_n}_0))]\nonumber\\
&=&\lim_{t\downarrow0}\frac{2}{t}E_{(f_{k,l}u_k)\cdot m}[(\widetilde{u_k}(X^{V_n}_t)-\widetilde{u_k}(X^{V_n}_0))(\tilde{w}(X^{V_n}_t)-\tilde{w}(X^{V_n}_0))]\nonumber\\
  &&+\lim_{t\downarrow0}\frac{1}{t}E_{f_{k,l}\cdot m}[(\widetilde{u_k}(X^{V_n}_t)-\widetilde{u_k}(X^{V_n}_0))^2(\tilde{w}(X^{V_n}_t)-\tilde{w}(X^{V_n}_0))]\nonumber\\
&:=&\lim_{t\downarrow0}[{I(t)+II(t)}].
\end{eqnarray}

By (\ref{aum1}), Theorem \ref{thm2.25}(iii) in the Appendix and (\ref{e4.27*}), we get
\begin{eqnarray}\label{m01}
\lim_{t\downarrow0}I(t)
&=&\lim_{t\downarrow0}\frac{2}{t}E_{(f_{k,l}u_k)\cdot m}(<M^{n,[u_k]},M^{n,[w]}>_t)\nonumber\\
&=&\lim_{t\downarrow0}\frac{2}{t}\int_{0}^t<\mu^{(n)}_{<u_k,w>},\widetilde{\hat{T}^{V_n}_s(f_{k,l}u_k)}>ds\nonumber\\
&=&\lim_{t\downarrow0}\frac{2}{t}\int_{0}^t[{\cal E}(u_k,w\hat{T}^{V_n}_s(f_{k,l}u_k))+{\cal E}(w,u_k\hat{T}^{V_n}_s(f_{k,l}u_k))\nonumber\\
& &\ \ \ \ \ \ \ \ \ \ -{\cal E}(u_kw,\hat{T}^{V_n}_s(f_{k,l}u_k))]ds.
\end{eqnarray}
By \cite[Theorem 3.4]{f}, $\hat{T}^{V_n}_s(f_{k,l}u_k)\rightarrow f_{k,l}u_k$ in $D(\mathcal{E})_{V_n}$ as $s\rightarrow0$. Furthermore, by Assumption \ref{assum1}, \cite[Corollary I.4.15]{MR92} and the fact that $|e^{-s}\hat{T}^{V_n}_s(f_{k,l}u_k)|\le k\|u_k\|_{\infty}\bar{h}_n$, $s>0$, we can show that
$\sup_{s>0}\mathcal{E}(w\hat{T}^{V_n}_s(f_{k,l}u_k),w\hat{T}^{V_n}_s(f_{k,l}u_k))<\infty$. Thus, we obtain by \cite[Lemma I.2.12]{MR92} that $w\hat{T}^{V_n}_s(f_{k,l}u_k)\rightarrow wf_{k,l}u_k$
weakly in $D(\mathcal{E})_{V_n}$ as $s\rightarrow0$. Similarly, we get $u_k\hat{T}^{V_n}_s(f_{k,l}u_k)\rightarrow u_kf_{k,l}u$
weakly in $D(\mathcal{E})_{V_n}$ as $s\rightarrow0$. Therefore, by (\ref{m01}) and (\ref{e4.27*}), we get
\begin{equation}\label{q01}
\lim_{t\downarrow0}I(t)=2\int_{V_n}\widetilde{f_{k,l}}\widetilde{u_k}d\mu^{(n)}_{<u_k,w>}.
\end{equation}

Note that
\begin{eqnarray}\label{q04}
II(t)
&=&\frac{1}{t}E_{f_{k,l}\cdot m}[(M^{n,[u_k],c}_t)^2M^{n,[w],c}_t]
   +\frac{1}{t}E_{f_{k,l}\cdot m}[(M^{n,[u_k],k}_t)^2M^{n,[w],k}_t]\nonumber\\
   &:=&III(t)+IV(t).
  \end{eqnarray}

By Burkholder-Davis-Gunday inequality, we get
\begin{eqnarray}\label{q02}
\lim_{t\downarrow0}III(t)
&\le&(\lim_{t\downarrow0}\frac{1}{t}E_{f_{k,l}\cdot m}[(M^{n,[u_k],c}_t)^4])^{1/2}(\lim_{t\downarrow0}\frac{1}{t}E_{f_{k,l}\cdot m}[<M^{n,[v],c}>_t])^{1/2}\nonumber\\
&\le&C(2ke^{V_n}(M^{n,[v]}))^{1/2}(\lim_{t\downarrow0}\frac{1}{t}E_{f_{k,l}\cdot m}[<M^{n,[u_k],c}>^2_t])^{1/2}
\end{eqnarray}
for some constant $C>0$, which is independent of $t$.

By Theorem \ref{thm2.25}(i) in the Appendix, for any $\delta>0$, we get
\begin{eqnarray}\label{n11}
& &\lim_{t\downarrow0}\frac{1}{t}E_{f_{k,l}\cdot m}[<M^{n,[u_k],c}>^2_t]\nonumber\\
&&\ \ \ \ =\lim_{t\downarrow0}\frac{2}{t}E_{f_{k,l}\cdot m}[\int_0^t<M^{n,[u_k],c}>_{(t-s)}\circ\theta_s d<M^{n,[u_k],c}>_s]\nonumber\\
&&\ \ \ \ =\lim_{t\downarrow0}\frac{2}{t}E_{f_{k,l}\cdot m}[\int_0^tE_{X^{V_n}_s}[<M^{n,[u_k],c}>_{(t-s)}]d<M^{n,[u_k],c}>_s]\nonumber\\
&&\ \ \ \ \le2<E_{\cdot}[<M^{n,[u_k]}>_{\delta}]\cdot\mu^{(n)}_{<u_k>},\widetilde{f_{k,l}}>.
\end{eqnarray}
Note that by our choice of $u_k$, there exists a constant $C_k>0$ such that $E_x(<M^{n,[u_k]}>_{\delta})=E_x[(M^{n,[u_k]}_{\delta})^2]=E_x[(\widetilde{u_k}(X^{V_n}_{\delta})-\widetilde{u_k}(X^{V_n}_0)-\int_0^{\delta}L^{V_n}u_k(X^{V_n}_s)ds)^2]\le C_k$ for any $\delta\le 1$ and ${\cal
E}{\textrm{-}q.e.}\ x\in V_n$. Letting $\delta\rightarrow 0$, by (\ref{n11}), the dominated convergence theorem and (\ref{q02}), we get
\begin{equation}\label{q111}\lim_{t\downarrow0}III(t)=0.
\end{equation}

By \cite[Theorem II.33, integration by parts (page 68) and Theorem II.28]{Pr},  we get
\begin{eqnarray}\label{q2222}
IV(t)
&=&\frac{1}{t}E_{f_{k,l}\cdot m}[I_{\zeta^{(n)}\le t}\{-(\widetilde{u_k}^2\tilde{w})(X^{V_n}_{\zeta^{(n)}-})
\nonumber\\
& &\ \ \ \ +2(\widetilde{u_k}\tilde{w})(X^{V_n}_{\zeta^{(n)}-})
(\widetilde{u_k}(X^{V_n}_{\zeta^{(n)}-})I_{\zeta^{(n)}\le t})^p\nonumber\\
& &\ \ \ \ +(\widetilde{u_k}^2)(X^{V_n}_{\zeta^{(n)}-})
(\tilde{w}(X^{V_n}_{\zeta^{(n)}-})I_{\zeta^{(n)}\le t})^p\}]\nonumber\\
&=&\frac{1}{t}E_{f_{k,l}\cdot m}[-((\widetilde{u_k}^2\tilde{w})(X^{V_n}_{\zeta^{(n)}-})I_{\zeta^{(n)}\le t})^p
\nonumber\\
& &\ \ \ \ +2((\widetilde{u_k}\tilde{w})(X^{V_n}_{\zeta^{(n)}-})I_{\zeta^{(n)}\le t})^p
M^{n,[u_k],k}_t\nonumber\\
& &\ \ \ \ +(\widetilde{u_k}^2(X^{V_n}_{\zeta^{(n)}-})I_{\zeta^{(n)}\le t})^p
M^{n,[w],k}_t]\nonumber\\
&\le&\frac{1}{t}E_{f_{k,l}\cdot m}[-((\widetilde{u_k}^2\tilde{w})(X^{V_n}_{\zeta^{(n)}-})I_{\zeta^{(n)}\le t})^p]
\nonumber\\
& &\ \ \ \ +\frac{2}{t}E^{1/2}_{f_{k,l}\cdot m}[\{((\widetilde{u_k}\tilde{w})(X^{V_n}_{\zeta^{(n)}-})I_{\zeta^{(n)}\le t})^p\}^2]
E^{1/2}_{f_{k,l}\cdot m}[<M^{n,[u_k],k}>_t]\nonumber\\
& &\ \ \ \ +E^{1/2}_{f_{k,l}\cdot m}[\{(\widetilde{u_k}^2(X^{V_n}_{\zeta^{(n)}-})I_{\zeta^{(n)}\le t})^p\}^2]
E^{1/2}_{f_{k,l}\cdot m}[<M^{n,[w],k}>_t].
\end{eqnarray}
By Theorem \ref{thm2.25}(i) in the Appendix, (\ref{l00})-(\ref{n02}), we obtain that for $\psi_1,\psi_2\in D(\mathcal{E})_{V_n,b}$,
\begin{eqnarray}\label{q22}
\lim_{t\downarrow0}\frac{1}{t}E_{f_{k,l}\cdot m}[((\widetilde{\psi_1}\widetilde{\psi_2})(X^{V_n}_{\zeta^{(n)}-})I_{\zeta^{(n)}\le t})^p]
&=&\int_{V_n}\widetilde{f_{k,l}}d\mu^{n,k}_{<\psi_1,\psi_2>}\nonumber\\
&=&\int_{V_n} \widetilde{f_{k,l}}\widetilde{\psi_1}\widetilde{\psi_2}dK^{(n)}
\end{eqnarray}
and
\begin{equation}\label{4.28}
\lim_{t\downarrow0}\frac{1}{t}E_{f_{k,l}\cdot m}[<M^{n,[\psi_1],k}>_t]
=\int_{V_n}\widetilde{f_{k,l}}d\mu^{n,k}_{<\psi_1>}.
\end{equation}
Furthermore, for any $\delta>0$,
\begin{eqnarray}\label{bn11}
& &\lim_{t\downarrow0}\frac{1}{t}E_{f_{k,l}\cdot m}[\{((\widetilde{\psi_1}\widetilde{\psi_2})(X^{V_n}_{\zeta^{(n)}-})I_{\zeta^{(n)}\le t})^p\}^2]\nonumber\\
&&\ \ \ \ =\lim_{t\downarrow0}\frac{2}{t}E_{f_{k,l}\cdot m}[\int_0^t((\widetilde{\psi_1}\widetilde{\psi_2})(X^{V_n}_{\zeta^{(n)}-})I_{\zeta^{(n)}\le (t-s)})^p\circ\theta_sd((\widetilde{\psi_1}\widetilde{\psi_2})(X^{V_n}_{\zeta^{(n)}-})I_{\zeta^{(n)}\le s})^p]\nonumber\\
&&\ \ \ \ =\lim_{t\downarrow0}\frac{2}{t}E_{f_{k,l}\cdot m}[\int_0^tE_{X_s^{V_n}}[
((\widetilde{\psi_1}\widetilde{\psi_2})(X^{V_n}_{\zeta^{(n)}-})I_{\zeta^{(n)}\le (t-s)})^p]d((\widetilde{\psi_1}\widetilde{\psi_2})(X^{V_n}_{\zeta^{(n)}-})I_{\zeta^{(n)}\le s})^p]\nonumber\\
&&\ \ \ \ \le<E_{\cdot}[(\widetilde{|\psi_1\psi_2|}(X^{V_n}_{\zeta^{(n)}-})I_{\zeta^{(n)}\le \delta})^p]\cdot\mu^{n,k}_{<|\psi_1|,|\psi_2|>},\widetilde{f_{k,l}}>\nonumber\\
&&\ \ \ \ =<E_{\cdot}[\widetilde{|\psi_1\psi_2|}(X^{V_n}_{\zeta^{(n)}-})I_{\zeta^{(n)}\le \delta}]\cdot\mu^{n,k}_{<|\psi_1|,|\psi_2|>},\widetilde{f_{k,l}}>.
\end{eqnarray}
Letting $\delta\rightarrow 0$, by (\ref{bn11}) and the dominated convergence theorem, we get
\begin{equation}\label{bn12}
\lim_{t\downarrow0}\frac{1}{t}E_{f_{k,l}\cdot m}[\{((\widetilde{\psi_1}\widetilde{\psi_2})(X^{V_n}_{\zeta^{(n)}-})I_{\zeta^{(n)}\le t})^p\}^2]=0.
\end{equation}

By (\ref{q2222})-(\ref{4.28}) and (\ref{bn12}), we get
\begin{equation}\label{bn13}\lim_{t\downarrow0}IV(t)=-\int_{V_n} \widetilde{f_{k,l}}\widetilde{u_k}^2\tilde{w}dK^{(n)}.
\end{equation}
Therefore, the proof is completed by (\ref{q00}), (\ref{q01}), (\ref{q04}), (\ref{q111}) and (\ref{bn13}).
\end{proof}

\begin{rem}\label{rem1}
When deriving formula (\ref{m00}) for non-symmetric Markov processes,
we cannot apply Theorem \ref{thm2.25}(vi) or (vii) in the Appendix of this paper to smooth measures which
are not of finite energy integral. To overcome that difficulty and obtain
(\ref{m00}) in the semi-Dirichlet forms setting, we have to make some extra
efforts as shown in the above proof. The proof uses some ideas of
\cite[Theorem 5.4]{Kim} and \cite[Theorem 5.3.2]{oshima}.
\end{rem}

\begin{thm}\label{thm4.6}
Let $m\in\mathbb{N}$, $\Phi\in C^1(\mathbb{R}^m)$ with $\Phi(0)=0$, and $u=(u_1,u_2,\dots, u_m)$ with $u_i\in D(\mathcal{E})_{V_n,b}$, $1\le i\le m$. Then
$\Phi(u)\in D(\mathcal{E})_{V_n,b}$
and for any $v\in D(\mathcal{E})_{V_n,b}$,
\begin{eqnarray}\label{e4.28}
d\mu^{n,c}_{<\Phi(u),v>}=\sum_{i=1}^m\Phi_{x_i}(\tilde{u})d\mu^{n,c}_{<u_i,v>}.
\end{eqnarray}
\end{thm}
\begin{proof}
$\Phi(u)\in D(\mathcal{E})_{V_n,b}$ is a direct consequence of Assumption \ref{assum1} and the corresponding property of Dirichlet form. Below we only prove (\ref{e4.28}). Let $v\in D(\mathcal{E})_{V_n,b}$. Then (\ref{e4.28}) is equivalent to
\begin{eqnarray}\label{new11.1}
\int_{V_n}\tilde{f}\bar{h}_nd\mu^{n,c}_{<\Phi(u),v>}=\sum_{i=1}^m\int_{V_n}\tilde{f}\bar{h}_n
\Phi_{x_i}(\tilde{u})d\mu^{n,c}_{<u_i,v>},\ \ \forall f\in D(\mathcal{E})_{V_n,b}.
\end{eqnarray}

Let $\mathcal{A}$ be the family of all $\Phi\in C^1(\mathbb{R}^m)$
satisfying (\ref{e4.28}). If $\Phi,\Psi\in \mathcal{A}$, then $\Phi\Psi\in \mathcal{A}$ by Theorem \ref{th4.5}. Hence $\mathcal{A}$ contains all polynomials vanishing at the origin. Let $O$ be a finite cube containing the range of $u(x)=(u_1(x),\dots,u_m(x))$. We take a sequence $\{\Phi^k\}$ of polynomials vanishing at the origin such that $\Phi^k\rightarrow\Phi$, $\Phi^k_{x_i}\rightarrow\Phi_{x_i}$, $1\leq i\leq m$, uniformly on $O$.
By Assumption \ref{assum1} and \cite[(3.2.27)]{Fu94},
$\Phi^k(u)$ converges to $\Phi(u)$ w.r.t. $\mathcal{E}^{V_n}_1$ as $k\rightarrow\infty$.
Then, by (\ref{e3.21}), we get
\begin{eqnarray*}
& &|\int_{V_n}\tilde{f}\bar{h}_nd\mu^{n,c}_{<\Phi(u),v>}
-\int_{V_n}\tilde{f}\bar{h}_nd\mu^{n,c}_{<\Phi^k(u),v>}|\nonumber\\
&&\ \ \ \ \ \ \leq\|f\|_{\infty}
|\int_{V_n}\bar{h}_nd\mu^{n,c}_{<\Phi(u)-\Phi^k(u)>}|^{1/2}
|\int_{V_n}\bar{h}_nd\mu^{n,c}_{<v>}|^{1
  /2}\nonumber\\
&&\ \ \ \ \ \ \leq\|f\|_{\infty}
|\int_{V_n}\bar{h}_nd\mu^{(n)}_{<\Phi(u)-\Phi^k(u)>}|^{1/2}
|\int_{V_n}\bar{h}_nd\mu^{(n)}_{<v>}|^{1
  /2}\nonumber\\
&&\ \ \ \ \ \ =2\|f\|_{\infty} e^{V_n}(M^{n,[\Phi(u)-\Phi^k(u)]})^{1/ 2}e^{V_n}(M^{n,[v]})^{1/2}\nonumber\\
&&\ \ \ \ \ \ \leq2\|f\|_{\infty}e^{V_n}(M^{n,[v]})^{1/ 2}[KC_n\mathcal{E}_{1}^{{V_n}}(\Phi(u)-\Phi^k(u),\Phi(u)-\Phi^k(u))^{1/2}\nonumber\\
&&\ \ \ \ \ \ \ \ \ \ \cdot(
     \|\Phi(u)-\Phi^k(u)\|_{\infty}\mathcal{E}_{1}^{{V_n}}(\bar{h}_{n},\bar{h}_{n})^{1/2}\\
     &&\ \ \ \ \ \ \ \ \ \ \ \ \ +\|\bar{h}_{n}\|_{\infty}\mathcal{E}_{1}^{{V_n}}(\Phi(u)-\Phi^k(u),\Phi(u)-\Phi^k(u))^{1/2})]
     ^{1/2}.\ \ \ \ \ \ \ \ \
\end{eqnarray*}
Hence
$$
\int_{V_n}\tilde{f}\bar{h}_nd\mu^{n,c}_{<\Phi(u),v>}=\lim_{k\rightarrow\infty}
\int_{V_n}\tilde{f}\bar{h}_nd\mu^{n,c}_{<\Phi^k(u),v>}.
$$
It is easy to see that
$$
\int_{V_n}\tilde{f}\bar{h}_n
\Phi_{x_i}(\tilde{u})d\mu^{n,c}_{<u_i,v>}=\lim_{k\rightarrow\infty}
\int_{V_n}\tilde{f}\bar{h}_n
\Phi^k_{x_i}(\tilde{u})d\mu^{n,c}_{<u_i,v>},\ \ 1\le i\le m.
$$
Therefore (\ref{new11.1}) holds.
\end{proof}
For $M,L\in\dot{\mathcal{M}}^{V_n}$, there exists a unique CAF $<M,L>$ of bounded variation such that
\begin{eqnarray*}
E_x(M_tL_t)=E_x(<M,L>_t),\ \ t\ge0,\ {\cal
E}{\textrm{-}q.e.}\ x\in V_n.
\end{eqnarray*}
Denote by $\mu^{(n)}_{<M,L>}$ the Revuz
 measure of $<M,L>$. Then, similar to
\cite[Lemma 5.6.1]{Fu94}, we can prove the following lemma.
\begin{lem}\label{lem4.7}
If $f\in L^{2}(V_n;\mu^{(n)}_{<M>})$ and $g\in L^{2}(V_n;\mu^{(n)}_{<L>})$, then $fg$ is integrable
w.r.t. $|\mu^{(n)}_{<M,L>}|$ and
\begin{eqnarray*}
(\int_{V_n}\mid fg\mid d\mid\mu^{(n)}_{<M,L>}\mid)^2\leq \int_{V_n}f^2 d\mu^{(n)}_{<M>}\int_{V_n}g^2 d\mu^{(n)}_{<L>}.
\end{eqnarray*}
\end{lem}
\begin{lem}\label{thm4.8}
Let $M\in\dot{\mathcal{M}}^{V_n}$ and $f\in L^2(V_n;\mu^{(n)}_{<M>})$. Then there exists a unique element
$f\cdot M\in\dot{\mathcal{M}}^{V_n}$ such that
$$
e^{V_n}(f\cdot M, L)=\frac{1}{2}\int_{V_n} f\bar{h}_nd\mu^{(n)}_{<M,L>}, \ \forall L\in\dot{\mathcal{M}}^{V_n}.
$$
The mapping $f\rightarrow f\cdot M$ is continuous and linear from
$L^2(V_n;\mu^{(n)}_{<M>})$ into the Hilbert space $(\dot{\mathcal{M}}^{V_n};e^{V_n})$.
\end{lem}
\begin{proof}
Let $L\in\dot{\mathcal{M}^{V_n}}$. Then, by Lemma \ref{lem4.7}, we get
\begin{eqnarray*}
\mid\frac{1}{2}\int_{V_n}f\bar{h}_nd\mu^{(n)}_{<M,L>}\mid
&\leq&\frac{1}{\sqrt{2}}(\int_{V_n} f^{2}\bar{h}_nd\mu^{(n)}_{<M>})^{1/2}
  ({1/2}\int_{V_n} \bar{h}_nd\mu^{(n)}_{<L>})^{1/2}\\
&\leq&\frac{\|\bar{h}_n\|_{\infty}}{\sqrt{2}}\parallel f\parallel_{L^{2}(V_n;\mu^{(n)}_{<M>})}\sqrt{e^{V_n}(L)}.
\end{eqnarray*}
Therefore, the proof is completed by Lemma \ref{lem3.4}.
\end{proof}

Similar to \cite[Lemma 5.6.2, Corollary 5.6.1 and Lemma 5.6.3]{Fu94}, we can prove the following two lemmas.
\begin{lem}\label{lem4.9}
Let $M,L\in\dot{\mathcal{M}}^{V_n}$. Then

(i) $d\mu^{(n)}_{<f\cdot M,L>}=fd\mu^{(n)}_{<M,L>}$ for $f\in L^{2}(V_n;\mu^{(n)}_{<M>})$.\\

(ii) $g\cdot(f\cdot M)=(gf)\cdot M$ for $f\in L^{2}(V_n;\mu^{(n)}_{<M>})$ and $g\in L^2(V_n;f^2d\mu^{(n)}_{<M>})$.\\

(iii) $e^{V_n}(f\cdot M,g\cdot L)={1\over 2}\int fg\bar{h}_nd\mu^{(n)}_{<M,L>}$ for $f\in L^{2}(V_n;\mu^{(n)}_{<M>})$ and $g\in L^{2}(V_n;\mu^{(n)}_{<L>})$.
\end{lem}\label{lem4.10*}
\begin{lem}\label{lemma4.7} The family
$\{\tilde{f}\cdot M^{u}\,|\,f\in D(\mathcal{E})_{V_n,b}\}$ is dense in $(\dot{\mathcal{M}}^{V_n},e^{V_n})$.
\end{lem}
\begin{thm}\label{3final}
Let $m\in\mathbb{N}$, $\Phi\in C^1(\mathbb{R}^m)$ with $\Phi(0)=0$, and $u=(u_1,u_2,\dots, u_m)$ with $u_i\in D(\mathcal{E})_{V_n,b}$, $1\le i\le m$. Then
\begin{eqnarray}\label{e4.32}
M^{[\Phi(u)],c}=\sum_{i=1}^m\Phi_{x_i}(u)\cdot M^{[u_i],c},\ \ P_{x}{\textrm{-}a.s.}\ \ {\rm for}\ {\cal
E}{\textrm{-}q.e.}\ x\in  V_n.
\end{eqnarray}
\end{thm}
\begin{proof} Let $v\in D(\mathcal{E})_{V_n,b}$ and $f,g\in D(\mathcal{E})_{V_n,b}$. Then, by Lemma \ref{lem4.9}(iii) and Theorem \ref{thm4.6}, we get
\begin{eqnarray*}
e^{V_n}(\tilde{f}\cdot M^{n,[\Phi(u)],c},\tilde{g}\cdot M^{n,[v]})
&=&{1\over 2}\int_{V_n} \tilde{f}\tilde{g}\bar{h}_nd\mu^{(n)}_{<M^{n,[\Phi(u)],c},M^{n,[v]}>}\\
&=&{1\over 2}\int_{V_n} \tilde{f}\tilde{g}\bar{h}_nd\mu^{n,c}_{<\Phi(u),v>}\\
&=&{1\over 2}\sum_{i=1}^m\int_{V_n} \tilde{f}\tilde{g}\bar{h}_n\Phi_{x_i}(u)d\mu^{n,c}_{<u_i,v>}\\
&=&{1\over 2}\sum_{i=1}^m\int_{V_n} \tilde{f}\tilde{g}\bar{h}_n\Phi_{x_i}(u)d\mu^{(n)}_{<M^{n,[u_i],c},M^{n,[v]}>}\\
&=&e^{V_n}(\sum_{i=1}^m(\tilde{f}\Phi_{x_i}(u))\cdot M^{n,[u_i],c},\tilde{g}\cdot M^{n,[v]}).
\end{eqnarray*}
By Lemma \ref{lemma4.7}, we get
$$
\tilde{f}\cdot M^{n,[\Phi(u)],c}=\sum_{i=1}^m(\tilde{f}\Phi_{x_i}(u))\cdot M^{n,[u_i],c},\ \ P_{x}{\textrm{-}a.s.}\ \ {\rm for}\ {\cal
E}{\textrm{-}q.e.}\ x\in  V_n.
$$
Therefore, (\ref{e4.32}) is satisfied by Lemma \ref{lem4.9}(ii), since $f\in D(\mathcal{E})_{V_n,b}$ is arbitrary.
\end{proof}

Let $M\in \dot{\mathcal{M}}_{loc}$. Then, there exist $\{V_n\},\{E_n\}\in\Theta$ and $\{M^n\,|\,M^n\in\dot{\mathcal{M}}^{V_n}\}$ such that $E_n\subset V_n$, $M_{t\wedge\tau_{E_n}}=M^{n}_{t\wedge\tau_{E_n}},\ t\ge0,\ n\in\mathbb{N}$. We define
$$
<M>_{t\wedge\tau_{E_n}}:=<M^n>_{t\wedge\tau_{E_n}};\ \ <M>_t:=\lim_{s\uparrow\zeta}<M>_s\ \ {\rm for}\ t\ge\zeta.
$$
Then, we can see that $<M>$ is well-defined and $<M>$ is a PCAF. Denote by $\mu_{<M>}$ the Revuz measure of $<M>$. We define
\begin{eqnarray*}
& &L^2_{loc}(E;\mu_{<M>})
 :=\{f\,|\, \exists\ \{V_n\},\{E_n\}\in\Theta\ {\rm and}\ \{M^n\,|\,M^n\in\dot{\mathcal{M}}^{V_n}\}\ \mbox{such that}\\
 &&\ \ \ \ \ \ \ E_n\subset V_n, M_{t\wedge\tau_{E_n}}=M^{n}_{t\wedge\tau_{E_n}},\ f\cdot I_{E_n}\in L^2(E_n;\mu^{(n)}_{<M^n>}),\ t\ge0,\ n\in\mathbb{N}\}
 \end{eqnarray*}
For $f\in L^2_{loc}(E;\mu_{<M>})$, we define $f\cdot M$ on $[\![0, \zeta[\![$ by
$$
(f\cdot M)_{t\wedge\tau_{E_n}}:=((f\cdot I_{E_n})\cdot M^n)_{t\wedge\tau_{E_n}},\ \ t\ge0,\ n\in\mathbb{N}.
$$
Then, we can see that $f\cdot M$ is well-defined and $f\cdot M\in {\mathcal{M}}^{[\![0, \zeta[\![}_{loc}$. Denote by $M^c$ the continuous part of $M$.

Finally, we obtain the main result of this section.
\begin{thm}\label{Ito} Suppose that
$(\mathcal{E},D(\mathcal{E}))$ is a quasi-regular local
semi-Dirichlet form on $L^{2}(E;m)$ satisfying Assumption \ref{assum1}.
Let $m\in\mathbb{N}$, $\Phi\in C^1(\mathbb{R}^m)$, and $u=(u_1,u_2,\dots, u_m)$ with $u_i\in D(\mathcal{E})_{loc}$, $1\le i\le m$. Then
$\Phi(u)\in D(\mathcal{E})_{loc}$ and
\begin{equation}\label{4final}
M^{[\Phi(u)],c}=\sum_{i=1}^m\Phi_{x_i}(u)\cdot M^{[u_i],c}\ {\rm on}\ [0,\zeta),\ \ P_{x}{\textrm{-}a.s.}\ \ {\rm for}\ {\cal
E}{\textrm{-}q.e.}\ x\in  E.
\end{equation}
\end{thm}
\begin{proof}
Since $1\in D(\mathcal{E})_{loc}$, $\Phi(u)\in D(\mathcal{E})_{loc}$ by Theorem \ref{thm4.6}. Hence (\ref{4final}) is a direct consequence of (\ref{e4.32}).
\end{proof}

\section[short
title]{Examples}\setcounter{equation}{0}\label{Sec:example}

In this section we investigate some concrete examples.

\begin{exa} We consider the following bilinear form
 \begin{eqnarray*}
\mathcal{E}(u,v)=\int^{1}_{0}u'v'dx+\int^{1}_{0}bu'vdx,\ \ u,v\in
D(\mathcal{E}):=H^{1,2}_{0}(0,1).
\end{eqnarray*}

\noindent (i) Suppose that $b(x)=x^{2}$. Then one can show that
$(\mathcal{E},D(\mathcal{E}))$ is a regular local semi-Dirichlet
form (but not a Dirichlet form) on $L^2((0,1);dx)$ (cf. \cite[Remark
2.2(ii)]{MR95}). Note that any $u\in D(\mathcal{E})$ is bounded and
$\frac{1}{2}-$H\"older continuous by the Sobolev embedding theorem.
Then we obtain Fukushima's decomposition,
$u(X_t)-u(X_0)=M^{[u]}_t+N^{[u]}_t$, by Lemma \ref{thm3.6}, where
$X$ is the diffusion process associated with
$(\mathcal{E},D(\mathcal{E}))$, $M^{[u]}$ is an MAF of finite energy
and $N^{[u]}$ is a CAF of zero energy.

\noindent (ii) Suppose that $b(x)=\sqrt{x}$. By \cite[Remark
2.2(ii)]{MR95}, $(\mathcal{E},D(\mathcal{E}))$ is a regular local
semi-Dirichlet form but not a Dirichlet form. Let $u\in
D(\mathcal{E})_{loc}$. Then we obtain Fukushima's decomposition
(\ref{new3}) by Theorem \ref{thm3.2}.

If $u\in D(\mathcal{E})$ satisfying ${\rm supp}[u]\subset (0,1)$,
then we may choose an open subset $V$ of $(0,1)$ such that ${\rm
supp}[u]\subset V\subset(0,1)$. Let $X^V$ be the part process of $X$
w.r.t. $V$. Then we obtain Fukushima's decomposition,
$u(X^V_t)-u(X^V_0)=M^{V,[u]}_t+N^{V,[u]}_t$, by Lemma \ref{thm3.6},
where $M^{V,[u]}$ is an MAF of finite energy and $N^{V,[u]}$ is a
CAF of zero energy w.r.t. $X^V$.
\end{exa}

\begin{exa}\label{e2} Let $d\geq3$, $U$ be an open subset of $\mathbb{R}^d$, $\sigma,\rho\in L^{1}_{loc}(U;dx)$, $\sigma,\rho>0$ $dx{\textrm{-}}a.e.$
For $u,v\in C^{\infty}_0(U)$, we define
\begin{eqnarray*}
\mathcal{E}_{\rho}(u,v)=\sum_{i,j=1}^d\int_U\frac{\partial
u}{\partial x_i}
 \frac{\partial v}{\partial x_j}\rho dx.
 \end{eqnarray*}
 Assume that
  $$(\mathcal{E}_{\rho},C^{\infty}_0(U)) \ \mbox{is closable on}\ L^2(U;\sigma dx).$$

Let $a_{ij},b_i,d_i\in L^1_{loc}(U;dx)$, $1\leq i,j\leq d$. For
$u,v\in C^{\infty}_{0}(U)$, we define
\begin{eqnarray*}
\mathcal{E}(u,v) &=&\sum_{i,j=1}^{d}\int_U\frac{\partial u}{\partial
     x_{i}}\frac{\partial u}{\partial
     x_{j}}a_{ij}dx+\sum_{i=1}^{d}\int_U
     \frac{\partial u}{\partial x_{i}}vb_{i}dx\\
     & &+\sum_{i=1}^{d}\int_U
     u\frac{\partial v}{\partial x_{i}}d_{i}dx
   +\int_U uvcdx.
\end{eqnarray*}
Set $\tilde{a}_{ij}:=\frac{1}{2}(a_{ij}+a_{ji})$,
$\check{a}_{ij}:=\frac{1}{2}(a_{ij}-a_{ji})$,
$\underline{b}:=(b_1,\dots, b_d)$, and $\underline{d}:=(d_1,\dots,
d_d)$. Define F to be the set of all functions $g\in
L^1_{loc}(U;dx)$ such that the distributional derivatives
$\frac{\partial g}{\partial x_i},\ 1\leq i\leq d$, are in
$L^1_{loc}(U;dx)$ such that $\|\nabla g\|(g\sigma)^{-\frac{1}{2}}\in
L^{\infty}(U;dx)$ or $\|\nabla
g\|^p(g^{p+1}\sigma^{p/q})^{-\frac{1}{2}}\in L^d(U;dx)$ for some
$p,q\in (1,\infty)$ with $\frac{1}{p}+\frac{1}{q}=1,\ p<\infty$,
where $\|\cdot\|$ denotes Euclidean distance in $\mathbb{R}^d$. We
say that a $\mathcal{B}(U)-$measurable function f has property
$(A_{\rho,\sigma})$ if one of the following conditions holds:

(i) $f(\rho\sigma)^{-\frac{1}{2}}\in L^{\infty}(U;dx)$.

(ii) $f^p(\rho^{p+1}\sigma^{p/q})^{-\frac{1}{2}}\in L^d(U,dx)$ for
some
 $p,q\in (1,\infty)$ with ${1\over p}+{1\over q}=1,\ p<\infty,$ and $\rho\in F$.

Suppose that

(C.I) There exists $\eta>0$ such that $ \sum_{i,j=1}^d
\tilde{a}_{ij}\xi_i\xi_j\ge\eta|\underline{\xi}|^2$, $\forall
\underline{\xi}=(\xi_1,\dots,\xi_d)\in\mathbb{R}^d$.

(C.II) $\check{a}_{ij}\rho^{-1}\in L^{\infty}(U;dx)$ for $1\le
i,j\le d$.

(C.III) For all $K\subset U$, $K$ compact,
$1_K\|\underline{b}+\underline{d}\|$
 and $1_Kc^{1/2}$ have property $(A_{\rho,\sigma}),$ and
 $(c+\alpha_0\sigma) dx-\sum_{i=1}^d\frac{\partial d_i}{\partial x_i}$
 is a positive measure on $\mathcal{B}(U)$ for some $\alpha_0\in (0,\infty)$.

(C.IV) $||\underline{b}-\underline{d}||$ has property
$(A_{\rho,\sigma})$.

(C.V) $\underline{b}=\underline{\beta}+\underline{\gamma}$ such that
$\|\underline\beta\|,
 \|\underline\gamma\|\in L^{1}_{loc}(U,dx)$,
 $(\alpha_0\sigma+c)dx-\sum_{1}^{d}{\partial \gamma_{i}\over \partial x_{i}}$
 is a positive measure on $\mathcal{B}(U)$ and
 $\|\underline{\beta}\|$
 has property $(A_{\rho,\sigma})$.

\noindent Then, by \cite[Theorem 1.2]{smuland}, there exists
$\alpha>0$ such that $(\mathcal{E}_{\alpha}, C^{\infty}_{0}(U))$ is
closable on $L^{2}(U;dx)$ and its closure
$(\mathcal{E}_{\alpha},D(\mathcal{E}_{\alpha}))$ is a regular local
semi-Dirichlet form on $L^{2}(U;dx)$. Define
$\eta_{\alpha}(u,u):=\mathcal{E}_{\alpha}(u,u)-\int\langle\triangledown
u,\underline{\beta}\rangle udx$
 for $u\in D(\mathcal{E}_{\alpha})$. By \cite[Theorem 1.2 (ii) and (1.28)]{smuland}, we know
$(\eta_{\alpha},D(\mathcal{E})_{\alpha})$ is a Dirichlet form and
there exists $C>1$ such that for any $u\in D(\mathcal{E}_{\alpha})$,
\begin{eqnarray*}
\frac{1}{C}\eta_{\alpha}(u,u)\leq \mathcal{E}_{\alpha}(u,u)\leq
C\eta_{\alpha}(u,u).
\end{eqnarray*}
 Let $X$ be the diffusion process associated with $(\mathcal{E}_{\alpha},D(\mathcal{E}_{\alpha}))$. Then, by Theorem \ref{thm3.2}, Fukushima's decomposition holds for any
$u\in D(\mathcal{E})_{loc}$. Moreover, the transformation formula (\ref{4final}) holds for local MAFs.
\end{exa}

\begin{exa}\label{e3} Let $S$ be a Polish space. Denote by ${\cal B}(S)$ the Borel $\sigma$-algebra of $S$. Let $E:={\cal M}_1(S)$ be the space of probability measures on $(S,{\cal B}(S))$.
For bounded ${\cal B}(S)$-measurable functions $f,g$ on $S$ and
$\mu\in E$, we define
$$
\mu(f):=\int_Sfd\mu,\ \ \langle
f,g\rangle_{\mu}:=\mu(fg)-\mu(f)\cdot\mu(g),\ \ \|f\|_{\mu}:=\langle
f,f\rangle_{\mu}^{1/2}.
$$
Denote by ${\cal F}C^{\infty}_b$ the family of all functions on $E$
with the following expression:
$$
u(\mu)=\varphi(\mu(f_1),\dots,\mu(f_k)),\ \ f_i\in C_b(S),1\le i\le
k,\varphi\in C^{\infty}_0(\mathbb{R}^k),k\in\mathbb{N}.
$$
Let $m$ be a finite positive measure on $(E,{\cal B}(E))$, where
${\cal B}(E)$ denotes the Borel $\sigma$-algebra of $E$. We suppose
that ${\rm supp}[m]=E$. Let $b:S\times E\rightarrow\mathbb{R}$ be a
measurable function such that
$$
\sup_{\mu\in E}\|b(\mu)\|_{\mu}<\infty,
$$
where $b(\mu)(x):=b(x,\mu)$.

For $u,v\in {\cal F}C^{\infty}_b$, we define
$$
{\cal E}^b(u,v):=\int_E(\langle \nabla u(\mu),\nabla
       v(\mu)\rangle_{\mu}+\langle b(\mu),\nabla u(\mu)\rangle_{\mu}v(\mu))m(d\mu),
$$
where
$$
\nabla u(\mu):=(\nabla_xu(\mu))_{x\in
S}:=\left(\left.\frac{d}{ds}u(\mu+s\varepsilon_x)\right|_{s=0}\right)_{x\in
S}.
$$
We suppose that $({\cal E}^0, {\cal F}C^{\infty}_b)$ is closable on
$L^2(E;m)$. Then, by \cite[Theorem 3.5]{ORS}, there exists
$\alpha>0$ such that $(\mathcal{E}^b_{\alpha}, {\cal
F}C^{\infty}_b)$ is closable on $L^{2}(E;m)$ and its closure
$(\mathcal{E}^b_{\alpha},D(\mathcal{E}^b_{\alpha}))$ is a
quasi-regular local semi-Dirichlet form on $L^{2}(E;m)$. Moreover,
by \cite[Lemma 2.5]{ORS}, there exists $C>1$ such that for any $u\in
D(\mathcal{E}^b_{\alpha})$,
\begin{eqnarray*}
\frac{1}{C}\mathcal{E}^0_{\alpha}(u,u)\leq
\mathcal{E}^b_{\alpha}(u,u)\leq C\mathcal{E}^0_{\alpha}(u,u).
\end{eqnarray*}
Let $X$ be the diffusion process associated with
$(\mathcal{E}^b_{\alpha},D(\mathcal{E}^b_{\alpha}))$, which is a
Fleming-Viot type process with interactive selection. Then, by
Theorem \ref{thm3.2}, Fukushima's decomposition holds for any $u\in
D(\mathcal{E}^b)_{loc}$. Moreover, the transformation formula (\ref{4final}) holds for local MAFs.
 \end{exa}

\section[short title]{Appendix: some results on potential theory and PCAFs
for semi-Dirichlet forms}\label{Appendix}
\setcounter{equation}{0}Let $E$ be a metrizable Lusin space and $m$ be a $\sigma$-finite positive
measure on its Borel $\sigma$-algebra ${\cal B}(E)$. Suppose that
$({\cal E},D({\cal E}))$ is a quasi-regular semi-Dirichlet form on
$L^2(E;m)$. Let $K>0$ be a continuity constant of $({\cal E},D({\cal
E}))$, i.e.,
\begin{equation}\label{K}
|{\cal E}_1(u,v)|\le K{\cal E}_1(u,u)^{1/2}{\cal E}_1(v,v)^{1/2},\ \
\forall u,v\in D({\cal E}).
\end{equation}
Denote by $(T_{t})_{t\geq0}$ and $(G_{\alpha})_{\alpha\geq0}$ (resp.
$(\hat{T}_{t})_{t\geq0}$ and $(\hat{G}_{\alpha})_{\alpha\geq0}$) the
semigroup and resolvent (resp. co-semigroup and co-resolvent)
associated with $({\cal E},D({\cal E}))$. Then there exists an
$m$-tight special standard process ${\bf M}=(\Omega,{\cal F},({\cal
F}_t)_{t\ge 0}, (X_t)_{t\ge 0},(P_x)_{x\in E_{\Delta}})$ which is
properly associated with $({\cal E},D({\cal E}))$
(cf. \cite[Theorem 3.8]{MR95}). It is known that any quasi-regular
semi-Dirichlet form is quasi-homeomorphic to a {regular}
semi-Dirichlet form (cf. \cite[Theorem 3.8]{HC06}). {By
quasi-homeomorphism and the transfer method (cf. \cite{CMR} and \cite[VI,
especially, Theorem VI.1.6]{MR92}), without loss of generality we
can restrict to Hunt processes when we discuss the AFs of ${\bf
M}$.}

Let $A\subset E$ and $f\in D(\mathcal{E})$. Denote by $f_{A}$ (resp.
$\hat{f}_A$) the 1-balayaged (resp. 1-cobalayaged) function of $f$
on $A$. \emph{We fix $\phi\in L^{2}(E;m)$
with
 $0<\phi\leq1$ $m{\textrm{-}}a.e.$ and set $h=G_{1}\phi$, $\hat{h}=\hat{G}_{1}\phi$.} Define for $U\subset
 E$, $U$ open,
 $${\rm cap}_{\phi}(U):=(h_{U},\phi)$$
and for any $A\subset E$,
 $${\rm cap}_{\phi}(A):=\inf\{{\rm cap}_{\phi}(U)\,|\,A\subset U, U\ {\rm open}\}.$$
 Hereafter, $(\cdot,\cdot)$ denotes the usual inner product of $L^2(E;m)$.
 By \cite[Theorem 2.20]{MR95}, we have
 $$
{\rm cap}_{\phi}(A)=(h_A,\phi)={\cal E}_1(h_A,\hat{G}_1\phi).
 $$

\begin{defin} A positive measure $\mu$ on $(E, {\cal B}(E))$ is said to be of {finite energy integral}, denoted by $S_0$, if $\mu(N)=0$ for each ${\cal E}$-exceptional set $N\in {\cal B}(E)$ and there exists a positive constant $C$ such that
$$
\int_E|\tilde{v}(x)|\mu(dx)\le C\mathcal{E}_{1}(v,v)^{1/2},\ \
\forall v\in {D}({\cal E}).
$$
\end{defin}

\begin{rem}\label{r1} (i)
Assume that $({\cal E},D({\cal E}))$ is a regular semi-Dirichlet
form. Let $\mu$ be a positive Radon measure on $E$ satisfying
$$
\int_E |v(x)|\mu(dx)\leq C{\cal E}_1(v,v)^{1/2},\ \ \forall v\in
C_0(E)\cap D(\cal{E})
$$
for some positive constant $C$, where $C_{0}(E)$ denotes the set of
all continuous functions on $E$ with compact supports. Then one can
show that $\mu$ charges no ${\cal E}$-exceptional set (cf.
\cite[Lemma 3.5]{HS10}) and thus $\mu\in S_0$.

(ii) Let $\mu\in S_0$ and $\alpha>0$. Then there exist unique
$U_{\alpha}\mu\in D(\mathcal{E})$ and  $\hat{U}_{\alpha}\mu\in
D(\mathcal{E})$ such that
\begin{equation}\label{e22}
\mathcal{E}_{\alpha}(U_{\alpha}\mu, v)=\int_E\tilde{v}(x)\mu(dx)=
\mathcal{E}_{\alpha}(v, \hat{U}_{\alpha}\mu).
\end{equation}
We call $U_{\alpha}\mu$ and $\hat{U}_{\alpha}\mu$ $\alpha$-potential
and $\alpha$-co-potential, respectively.

Let $u\in D({\cal E})$. By quasi-homeomorphism and similar to
\cite[Theorem 2.2.1]{Fu94} (cf. \cite[Lemma 1.2]{HS10}), one can
show that the following conditions are equivalent to each other:

(i) $u$ is $\alpha$-excessive (resp. $\alpha$-co-excessive).

(ii) $u$ is an $\alpha$-potential (resp. $\alpha$-co-potential).

 (iii) ${\cal E}_{\alpha}(u,v)\geq 0$ (resp. ${\cal E}_{\alpha}(v,u)\geq 0$), $\forall v\in D({\cal E}),\ v\geq 0$.
\end{rem}

\begin{thm}\label{thm34} Define
$$
\hat{S}^{*}_{00}:=\{\mu\in S_{0}\,|\,\hat{U}_{1}\mu\leq
c\hat{G}_{1}\phi\ \mbox{for some constant}\ \ c>0\}.
$$
Let $A\in {\cal B}(E)$. If $\mu(A)=0$ for all $\mu\in
\hat{S}^{*}_{00}$, then cap$_{\phi}(A)=0$.
\end{thm}
\begin{proof} By quasi-homeomorphism, without loss of generality, we suppose that $({\cal E}, D({\cal E}))$ is a regular semi-Dirichlet form.
Assume that $A\in {\cal B}(E)$ satisfying $\mu(A)=0$ for all $\mu\in
\hat{S}^{*}_{00}$. We will prove that cap$_{\phi}(A)=0$.

\noindent {\it Step 1.} We first show that $\mu(A)=0$ for all
$\mu\in S_0$. Suppose that $\mu\in S_0$. By \cite[Proposition
4.13]{MR1995}, there exists an ${\cal E}$-nest $\{F_k\}$ of compact subsets of $E$ such that
$\widetilde{\hat{G}_{1}\phi},\ \widetilde{\hat{U}_{1}\mu}\in
C(\{F_k\})$ and $\widetilde{\hat{G}_{1}\phi}>0$ on $F_k$ for each
$k\in\mathbb{N}$. Then, there exists a sequences of positive
constants $\{a_k\}$ such that
$$
\widetilde{\hat{U}_{1}\mu}\le a_k \widetilde{\hat{G}_{1}\phi}\ \
{\rm on}\ F_k\ {\rm for\ each}\ k\in\mathbb{N}.
$$

Define $u_k=\hat{U}_1(I_{F_k}\cdot\mu)$ and set $v_k=u_k\wedge a_k
{\hat{G}_{1}\phi}$ for $k\in\mathbb{N}$. Then $\widetilde{u_k}\le
\widetilde{\hat{U}_1\mu}\le a_k\widetilde{\hat{G}_{1}\phi}$ ${\cal
E}{\textrm{-}}q.e.$ on $F_k$. By (\ref{e22}), we get
$$
{\cal
E}_1(v_k,u_k)=\int_{F_k}\widetilde{v_k}(x)\mu(dx)=\int_{F_k}\widetilde{u_k}(x)\mu(dx)={\cal
E}_1(u_k,u_k).
$$
Since $v_k$ is a 1-co-potential and $v_k\le u_k$
$m{\textrm{-}}a.e.$, ${\cal E}_1(v_k-u_k,v_k-u_k)={\cal
E}_1(v_k-u_k, v_k)-{\cal E}_1(v_k-u_k,u_k)\le 0$, proving that
$u_k=v_k\le a_k {\hat{G}_{1}\phi}$ $m{\textrm{-}}a.e.$ Hence
$I_{F_k}\cdot\mu\in \hat{S}^{*}_{00}$. Therefore $\mu(A)=0$ by the
assumption that $A$ is not charged by any measure in
$\hat{S}^{*}_{00}$.

\noindent {\it Step 2.} Suppose that cap$_{\phi}(A)>0$. By
\cite[Corollary 2.22]{MR95}, there exists a compact set $K\subset B$
such that cap$_{\phi}(K)>0$. Note that
$(\widehat{\hat{G}_{1}\phi})_{K}\in D(\mathcal{E})$ is
1-co-excessive. By Remark \ref{r1}(ii), there exists
$\mu_{(\widehat{\hat{G}_{1}\phi})_{K}}\in S_{0}$ such that
\begin{eqnarray}\label{e23}
{\rm
cap}_{\phi}(K)&=&\mathcal{E}_{1}((G_{1}\phi)_K,{\hat{G}_{1}\phi})\nonumber\\
&=&\mathcal{E}_{1}(G_{1}\phi,(\widehat{\hat{G}_{1}\phi})_{K})\nonumber\\
&=&\int_{E}\widetilde{G_{1}\phi}
d\mu_{(\widehat{\hat{G}_{1}\phi})_{K}}\nonumber\\
&\leq&\mu_{(\widehat{\hat{G}_{1}\phi})_{K}}(E).
\end{eqnarray}
For any $v\in C_{0}(K^{c})\cap D(\mathcal{E})$, we have $\int
\tilde{v}d\mu_{(\widehat{\hat{G}_{1}\phi})_{K}}=\mathcal{E}_{1}(v,(\widehat{\hat{G}_{1}\phi})_{K})=0$.
Since $C_{0}(K^{c})\cap D(\mathcal{E})$ is dense in
$C_{0}(K^{c})$, the support of $\mu_{\widehat{\hat{G}_{1}\phi}}$ is
contained in $K$. Thus, by (\ref{e23}), we get
$\mu_{\widehat{\hat{G}_{1}\phi}}(K)>0$. Therefore cap$_{\phi}(A)=0$
by Step 1.
\end{proof}

\begin{thm}\label{thm2.13}
The following conditions are equivalent for a positive measure $\mu$
on $(E, {\cal B}(E))$.

(i) $\mu\in S$.

 (ii) There exists an ${\cal E}$-nest $\{F_k\}$
satisfying $I_{F_{k}}\cdot\mu\in S_{0}$ for each $k\in \mathbb{N}$.
\end{thm}
\begin{proof} (ii) $\Rightarrow$ (i) is clear. We only prove (i) $\Rightarrow$
(ii). Let $(\tilde{\cal E}, D({\cal E}))$ be the symmetric part of
$({\cal E}, D({\cal E}))$. Then $(\tilde{\cal E}, D({\cal E}))$ is a
symmetric positivity preserving form. Denote by
$(\tilde{G}_{\alpha})_{\alpha\geq0}$ the resolvent associated with
$(\tilde{\cal E},D({\cal E}))$ and set
$\bar{h}:=\tilde{G}_1\varphi$. Then $(\tilde{\cal E}^{\bar{h}}_1,
D({\cal E}^{\bar{h}}))$ is a quasi-regular symmetric Dirichlet form
on $L^2(E;{\bar{h}}^2\cdot m)$ (the ${\bar{h}}$-transform of
$(\tilde{\cal E}_1, D({\cal E}))$).

By \cite[pages 838-839]{kuwae}, for an increasing sequence $\{F_k\}$
of closed sets, $\{F_k\}$ is an ${\cal E}$-nest if and only if it is
an $\tilde{\cal E}^{\bar{h}}_1$-nest. We select a compact
$\tilde{\cal E}^{\bar{h}}_1$-nest $\{F_k\}$ such that
$\widetilde{\bar{h}}$ is bounded on each $F_k$. Let $\mu\in S({\cal
E})$, the family of smooth measures w.r.t. $(E, {\cal B}(E))$. Then
$\mu\in S(\tilde{\cal E}^{\bar{h}}_1)$, the family of smooth
measures w.r.t. $(\tilde{\cal E}^{\bar{h}}_1, D({\cal
E}^{\bar{h}}))$. By \cite[Theorem 2.2.4]{Fu94} and
quasi-homeomorphism, we know that there exists a compact
$\tilde{\cal E}^{\bar{h}}_1$-nest (hence ${\cal E}$-nest) $\{J_k\}$
such that $I_{J_k}\cdot\mu\in S_0(\tilde{\cal E}^{\bar{h}}_1)$.
Then, there exists a sequence of positive constants $\{C_k\}$ such
that
$$
\int_E|\tilde{g}|I_{J_k}d\mu\le C_k{\tilde{\cal
E}^{\bar{h}}_1(g,g)}^{1/2}, \ \ \forall g\in D({\cal E}^{\bar{h}}).
$$

We now show that each $I_{F_k\cap J_k}\cdot\mu\in S_0({\cal E})$,
which will complete the proof. In fact, let $f\in D({\cal E})$. Then
$ \frac{f}{{\bar{h}}}\in D({\cal E}^{\bar{h}})$ and
\begin{eqnarray*}
\int_E|\tilde{f}|I_{F_k\cap J_k}d\mu&\le&
\|\,{\bar{h}}|_{F_k}\|_{\infty}\int_E|\frac{\tilde{f}}{\bar{h}}|I_{F_k\cap J_k}d\mu\\
&\le&
\|\,{\bar{h}}|_{F_k}\|_{\infty}\int_E|\frac{\tilde{f}}{\bar{h}}|I_{J_k}d\mu\\
&\le& \|\,\bar{h}|_{F_k}\|_{\infty}C_k{\tilde{\cal E}^{\bar{h}}_1(f/{\bar{h}},f/{\bar{h}})}^{1/2}\\
&=&\|\,{\bar{h}}|_{F_k}\|_{\infty}C_k{{\cal E}_1(f,f)}^{1/2}.
\end{eqnarray*}
Since $f\in D({\cal E})$ is arbitrary, this implies that $I_{F_k\cap
J_k}\cdot\mu\in S_0({\cal E})$.
\end{proof}

\begin{lem}\label{lem2.19} For any $u\in D(\mathcal{E})$, $\nu\in S_{0}$,
$0<T<\infty$ and $\varepsilon>0$,
\begin{eqnarray*}
P_{\nu}(\sup_{0\leq t\leq T}|\tilde{u}(X_{t})|>\varepsilon)
\leq\frac{2K^{5/2}e^{T}}{\varepsilon}{\mathcal{E}_{1}(u,u)}^{1/2}{\mathcal{E}_{1}(\hat{U}_{1}\nu,
\hat{U}_{1}\nu)}^{1/2}.
\end{eqnarray*}
\end{lem}
\begin{proof} We take an ${\cal E}$-quasi-continuous Borel version $\tilde{u}$ of $u$. Let $A=\{x\in E\,|\,
|\tilde{u}(x)|>\varepsilon\}$ and
$\sigma_{A}:=\inf\{t>0\,|\,X_{t}\in A\}$. By \cite[Theorem
4.4]{kuwae},
$H_{A}^{1}|u|:=E_{\cdot}[e^{-\sigma_{A}}|u|(X_{\sigma_{A}})]$ is an
${\cal E}$-quasi-continuous version of $|u|_{A}$. Then, by
\cite[Proposition 2.8(i) and (2.1)]{MR95}, we get
\begin{eqnarray*}
P_{\nu}(\sup_{0\leq t\leq T}|\tilde{u}(X_{t})|>\varepsilon)&
\leq&\frac{e^{T}E_{\nu}[e^{-\sigma_{A}}|u|(X_{\sigma_{A}})]}{\varepsilon}\\
&=&\frac{e^{T}}{\varepsilon}\int_E|u|_{A}d\nu\\
&=&\frac{e^{T}}{\varepsilon}\mathcal{E}_{1}(|u|_{A},\hat{U_{1}}\nu)\\
&\leq&\frac{Ke^{T}}{\varepsilon}\mathcal{E}_{1}(|u|_{A},|u|_{A})^{1/2}
     \mathcal{E}_{1}(\hat{U_{1}}\nu,\hat{U_{1}}\nu)^{1/2}\\
&\leq&\frac{K^{2}e^{T}}{\varepsilon}\mathcal{E}_{1}(|u|,|u|)^{1/2}
     \mathcal{E}_{1}(\hat{U_{1}}\nu,\hat{U_{1}}\nu)^{1/2}\\
&\leq&\frac{2K^{3}e^{T}}{\varepsilon}\mathcal{E}_{1}(u,u)^{1/2}
     \mathcal{E}_{1}(\hat{U_{1}}\nu,\hat{U_{1}}\nu)^{1/2}.
\end{eqnarray*}
\end{proof}

By Lemma \ref{lem2.19} and Theorem \ref{thm34}, similar to
\cite[Lemma 5.1.2]{Fu94}, we can prove the following lemma.

\begin{lem}\label{l343} Let $\{u_{n}\}$ be a sequence of ${\cal E}$-quasi continuous functions in
$D(\mathcal{E})$. If $\{u_{n}\}$ is an $\mathcal{E}_{1}$-Cauchy
sequence, then there exists a subsequence $\{u_{n_{k}}\}$ satisfying
the condition that for ${\cal E}{\textrm{-}}q.e.$ $x\in E$,
$$
P_{x}(u_{n_{k}}(X_{t}) \ \mbox{ converges uniformly in t on each
compact interval of}\ \  [0,\infty))=1.
$$
\end{lem}

In \cite{F01}, Fitzsimmons extended the smooth measure
characterization of PCAFs from the Dirichlet forms setting to the
semi-Dirichlet forms setting (see \cite[Theorem 4.22]{F01}). In
particular, the following proposition holds.

\begin{prop}\label{p23} (cf. \cite[Proposition 4.12]{F01}) For any
$\mu\in S_{0}$, there is a unique finite PCAF $A$ such that
$E_{x}(\int_{0}^{\infty}e^{-t}dA_{t})$ is an ${\cal
E}$-quasi-continuous version of $U_{1}\mu$.
\end{prop}

By Proposition \ref{p23} and Theorem \ref{thm2.13}, following the
arguments of \cite[Theorems 5.1.3 and 5.1.4]{Fu94} (with slight
modifications by virtue of \cite{MR95,MR1995,kuwae} and
\cite[Theorem 3.4]{f}), we can obtain the following theorem.

 \begin{thm}\label{thm2.25} Let $\mu\in S$ and $A$ be a PCAF. Then
 the following conditions are equivalent to each other:

(i) For any $\gamma$-co-excessive function $g$ $(\gamma\geq0)$ in
$D({\cal E})$ and $f\in \mathcal{B}^{+}(E)$,
\begin{equation}\label{e1.2}
 \lim_{t\downarrow0}\frac{1}{t}E_{g\cdot
 m}((fA)_{t})=<f\cdot\mu,\tilde{g}>.
\end{equation}

(ii)For any $\gamma$-co-excessive function $g$ $(\gamma\geq0)$  in
$D({\cal E})$ and $f\in \mathcal{B}^{+}(E)$,
$$
\alpha(g,U^{\alpha+\gamma}_{A}f)\ \uparrow\ <f\cdot\mu,\tilde{g}>,\
\ \alpha\uparrow\infty,
$$
where $U^{\alpha}_Af(x):=E_x(\int_0^{\infty}e^{-\alpha
t}f(X_t)dA_t)$.

 (iii) For any $t>0$, $g\in \mathcal{B}^{+}(E)\cap L^2(E;m)$ and $f\in \mathcal{B}^{+}(E)$,
 $$E_{g\cdot m}((fA)_{t})=\int^{t}_{0}<f\cdot\mu,\widetilde{\hat{T}_{s}g}>ds.$$

 (iv) For any $\alpha>0$, $g\in \mathcal{B}^{+}(E)\cap L^2(E;m)$ and $f\in \mathcal{B}^{+}(E)$,
 $$(g,U^{\alpha}_{A}f)=<f\cdot\mu,\widetilde{\hat{G}_{\alpha}g}>.$$

 When $\mu\in S_{0}$, each of the above four conditions is also
 equivalent to each of the following three conditions:

 (v) $U^{1}_{A}1$ is an ${\cal E}$-quasi-continuous version of $U_{1}\mu$.

(vi) For any $g\in\mathcal{B}^{+}(E)\cap D(\mathcal{E})$ and
$f\in\mathcal{B}^{+}_{b}(E)$,
$$
\lim_{t\downarrow0}\frac{1}{t}E_{g\cdot
m}((fA)_{t})=<f\cdot\mu,\tilde{g}>.
$$

(vii) For any $g\in\mathcal{B}^{+}(E)\cap D(\mathcal{E})$ and
$f\in\mathcal{B}^{+}_{b}(E)$,
$$
\lim_{\alpha\rightarrow\infty}\alpha(g,U^{\alpha}_{A}f)=<f\cdot\mu,\tilde{g}>.$$

The family of all equivalent classes of PCAFs and the family $S$ are
in one to one correspondence under the Revuz correspondence
(\ref{e1.2}).
\end{thm}

Given a PCAF $A$, we denote by $\mu_A$ the Revuz measure of $A$.

\begin{lem}\label{lnm} Let $A$ be a PCAF and $\nu\in
\hat{S}^{*}_{00}$. Then there exists a positive constant $C_{\nu}$
such that for any $t>0$,
$$
E_{\nu}(A_{t})\leq C_{\nu} (1+t)\int_E\widetilde{\hat{h}}d\mu_A.
$$
\end{lem}
\begin{proof} By Theorem \ref{thm2.13}, we may assume without loss of
generality that $\mu_A\in S_{0}$. Set $c_{t}(x)=E_{x}(A_{t})$.
Similar to \cite[page 137]{oshima}, we can show that $c_t\in
D(\mathcal{E})$ and for any $v\in D(\mathcal{E})$
$$\mathcal{E}(c_{t},v)=<\mu_A,v-\hat{T}_{t}v>.$$

Let $\nu\in \hat{S}^{*}_{00}$. Then
\begin{eqnarray*}
E_{\nu}(A_{t})
&=&<\nu, c_{t}>\\
&=&\mathcal{E}_{1}(c_{t},\hat{U}_{1}\nu)\\
&\leq&<\mu_A,\hat{U}_{1}\nu>+<c_{t},\hat{U}_{1}\nu>\\
&\leq&C_{\nu}[<\mu_A,\hat{h}>+E_{\hat{h}\cdot m}(A_{t})]
\end{eqnarray*}
for some constant $C_{\nu}>0$. Therefore the proof is completed by
(\ref{e1.2}).
\end{proof}

%%%%%%%%%%%%%%%%%%%%%%%%%%%%%%%%%%%%%%%%%%%%%%%%%%%%%%%%%%%%%%%%%%%%%%%%%%%%%%%%%%%%%%%%%%%%%%%%%%%%%%%%
\bigskip

{ \noindent {\bf\large Acknowledgments} \vskip 0.1cm  \noindent We are
grateful to the support of NSFC (Grant No. 10961012), 973 Project, Key Lab of CAS (Grant No. 2008DP173182), and
NSERC (Grant No. 311945-2008).}
%%%%%%%%%%%%%%%%%%%%%%%%%%%%%%%%%%%%%%%%%%%%%%%%%%%%%%%%%%%%%%%%%%%%%%%%%%%%%%%%%%%%%%%%%%%%%%%%%%%%%%

\end{document}